\documentclass{article}

\usepackage{geometry}
\usepackage{authblk}
\usepackage{amsmath}
\usepackage{graphicx}
\usepackage{algorithmic}
\usepackage[title]{appendix}
\usepackage{algorithm}
\usepackage{lmodern}
\usepackage{float}
\usepackage{subfigure}
\usepackage{booktabs}
\usepackage{appendix}
\usepackage{cite}
\usepackage[colorlinks=true,  allcolors=blue]{hyperref}

\begin{document}

\title{An Unsupervised Network Architecture Search Method for Solving Partial Differential Equations}

\author[1, \thanks{lq\_2019@mail.ustc.edu.cn}]{Qing Li}
\author[2, 3, \thanks{Corresponding author: jingrunchen@ustc.edu.cn}]{Jingrun Chen}
\affil[1]{School of Artificial Intelligence and Data Science, University of Science and Technology of China, China}    
\affil[2]{School of Mathematical Sciences and Suzhou Institute for Advanced Research, University of Science and Technology of China, China}
\affil[3]{Suzhou Big Data \& AI Research and Engineering Center, China}

\date{}

\maketitle

\begin{abstract}
Solving partial differential equations (PDEs) has been indispensable in scientific and engineering applications.  
Recently,  deep learning methods have been widely used to solve high-dimensional problems. The interplay between deep learning and numerical PDEs has resulted in a variety of new methods, one of which is the physics-informed neural network (PINN). Typically, a deep learning method has three main components: a neural network,  a loss function,  and an optimizer. While the construction of the loss function is rooted in the definition of solution space, how to choose a \textit{optimal} neural network is somewhat ad hoc,  leaving much room for improvement. In the framework of PINN,  we propose an unsupervised network architecture search method for solving PDEs, termed PINN-DARTS,  which applies the differentiable architecture search (DARTS) to find the optimal network architecture structure in a given set of neural networks. 
In this set,  the number of layers and the number of neurons in each layer can change. In the searching phase,  both network and architecture parameters are updated simultaneously, so the running time is close to that of PINN with a pre-determined network structure. Unlike available works,  our approach is unsupervised and purely based on the PDE residual without any prior usage of solutions. PINN-DARTS outputs the optimal network structure as well as the associated numerical solution. The performance of PINN-DARTS is verified on several benchmark PDEs,  including elliptic, parabolic, wave, and Burgers' equations. Compared to traditional architecture search methods, PINN-DARTS achieves significantly higher architectural accuracy. Another interesting observation is that both the solution complexity 
and the PDE type have a prominent impact on the optimal network architecture. Our study suggests that architectures with uneven widths from layer to layer may have superior performance across different solution complexities and different PDE types. 
\end{abstract}

\section{Introduction}
\indent \indent
Partial differential equations (PDEs) are widely used in scientific and engineering problems. In most scenarios, one has to look for numerical solutions due to the complexity of models in real applications. Popular methods, such as Finite Difference Method \cite{FDM}, Finite Element Method \cite{FEM}, and Finite Volume Method \cite{FVM}, have achieved significant success. Industrial software based on these methods provides a convenient yet powerful tool for academia and industry. However, when the dimensionality of the problem increases, the computational cost of these methods grows exponentially, leading to the curse of dimensionality.
\par
In recent years, with the rapid development of deep learning, the idea of using neural networks to solve PDEs has gained widespread attention. One such method is the physics-informed neural network (PINN) \cite{PINN}, which incorporates physical constraints into the loss function of neural networks to solve high-dimensional PDEs. Due to their simplicity and good experimental performance, PINN have become a representative method in the field of scientific machine learning \cite{1-1}, and are widely used to solve various problems, such as fluid mechanics \cite{1-2}, structural mechanics \cite{1-3}, and industrial design \cite{1-4}.
\par
Typically, the network architecture approximating the PDE solution is a feedforward neural network (FNN) \cite{PINN,2-6}. While there has been considerable research on loss functions and network errors for different problems using PINN \cite{2-1,2-2,2-3}, studies on the optimal network architecture are limited \cite{PINN,2-4,2-5}. Hyperparameters of the network architecture are pre-selected, including network width, network depth, activation function, and learning rate. Along another line, in the field of image processing, there are significant efforts for searching the optimal architecture, i.e., neural architecture search (NAS) \cite{3-1,3-6}. NAS searches for an optimal neural network architecture within a prescribed search space in an auto-design way. Successful approaches in NAS, include reinforcement learning \cite{3-2,3-7}, gradient-based optimization \cite{DARTS,3-3,3-8}, evolutionary algorithms \cite{3-4,3-9}, and Bayesian optimization \cite{3-5,3-10}. As a gradient-based method, differentiable architecture search (DARTS) \cite{DARTS} is a popular method in which both network parameters and architecture parameters are simultaneously optimized using gradient descent.

\par 
There are few attempts to apply the idea of NAS to search for the optimal architecture in PINN. For example, Auto-PINN \cite{Auto-PINN} uses a step-by-step decoupling strategy to search, optimizing one hyperparameter at each step with the others fixed to one or a few sets of options. 
Another example is the NAS-PINN \cite{NAS-PINN}, which expands the search space by adding masks to allow uneven layer widths for each layer of the FNN. The architecture searching phase is formulated as a two-layer optimization problem: the inner layer optimizes PINN network parameters, and the outer loop updates the hyperparameters of PINN. However, the metric used in NAS-PINN is the error between the numerical solution of PINN and the exact solution.
 
\par
In this paper,  we propose PINN-DARTS,  a class of unsupervised search methods for automating the design of optimal architectures in PINN. Thanks to the differentiability of DARTS, we optimize the architecture parameters and network parameters simultaneously using the gradient descent method in the searching phase. Hyperparameters to be optimized include the number of layers and the number of neurons in each layer. Unlike available works, PINN-DARTS is fully unsupervised, applicable to more practical problems. Compared to traditional search methods, PINN-DARTS performs better in search time and architectural accuracy.  We validate the superiority of PINN-DARTS across various PDEs with different solution complexities.

\par
The paper is organized as follows. Section \ref{section2} introduces the PINN and network search methods, including traditional search methods and DARTS. Our approach, PINN-DARTS, as well as methodological details, are described in Section \ref{section3}.
In Section \ref{section4}, numerical experiments are conducted on different PDEs to demonstrate the effectiveness and advantages of PINN-DARTS over other search methods. Conclusions and future work are drawn in Section \ref{section5}.

\section{Preliminaries}
\label{section2}

\subsection{Physics-Informed Neural Network}
\indent \indent

As a general framework for solving PDEs, physics-informed neural network (PINN) \cite{PINN} integrates the equation, the boundary condition, the initial condition, and the observation data (if available) into the loss function. 
Consider the following PDE problem
\begin{equation}
    \label{2.1}
    \left\{
        \begin{aligned}
        \mathcal{L}(u(\boldsymbol{x},t)) & = 0,  \quad & \boldsymbol{x} \in \Omega,  \; t \in [0,T],  \\
         \mathcal{B}(u(\boldsymbol{x},t)) & = b(\boldsymbol{x},t),  \quad & \boldsymbol{x} \in \partial \Omega,  \; t \in [0,T],  \\
         u(\boldsymbol{x},0) & = u_0(\boldsymbol{x}),  \quad & \boldsymbol{x} \in \Omega. 
        \end{aligned}
    \right. 
\end{equation}
The loss function in PINN is defined as the sum of the mean squared error (MSE) of residuals associated with the equation, the boundary condition, and the initial condition.
\begin{equation}
    \label{eqn:loss}
    {L}_{\text{PINN}} = \lambda_f {L}_{\text{PDE}} + \lambda_b {L}_{\text{BC}} + \lambda_i {L}_{\text{IC}}, 
\end{equation}
where
\begin{align*}
    {L}_{\text{PDE}} & = \frac{1}{N_f} \sum_{i=1}^{N_f} \left| \mathcal{L}(u(\boldsymbol{x}_i,t_i)) \right|^2, \\
    {L}_{\text{BC}} & = \frac{1}{N_b} \sum_{i=1}^{N_b} \left| u(\boldsymbol{x}_i,t_i) - b(\boldsymbol{x}_i,t_i) \right|^2, \\
    {L}_{\text{IC}} & = \frac{1}{N_i} \sum_{i=1}^{N_i} \left| u(\boldsymbol{x}_i,0) - u_0(\boldsymbol{x}_i) \right|^2.
\end{align*}
$\lambda_f$, $\lambda_b$, and $\lambda_b $ are the hyperparameters that can be tuned. Uniform sampling or random sampling is employed to generate sampling points $\{(\boldsymbol{x}_i, t_i)\}$.
\par
\begin{figure}[htbp]
    \centering
    \includegraphics[width=0.75\linewidth]{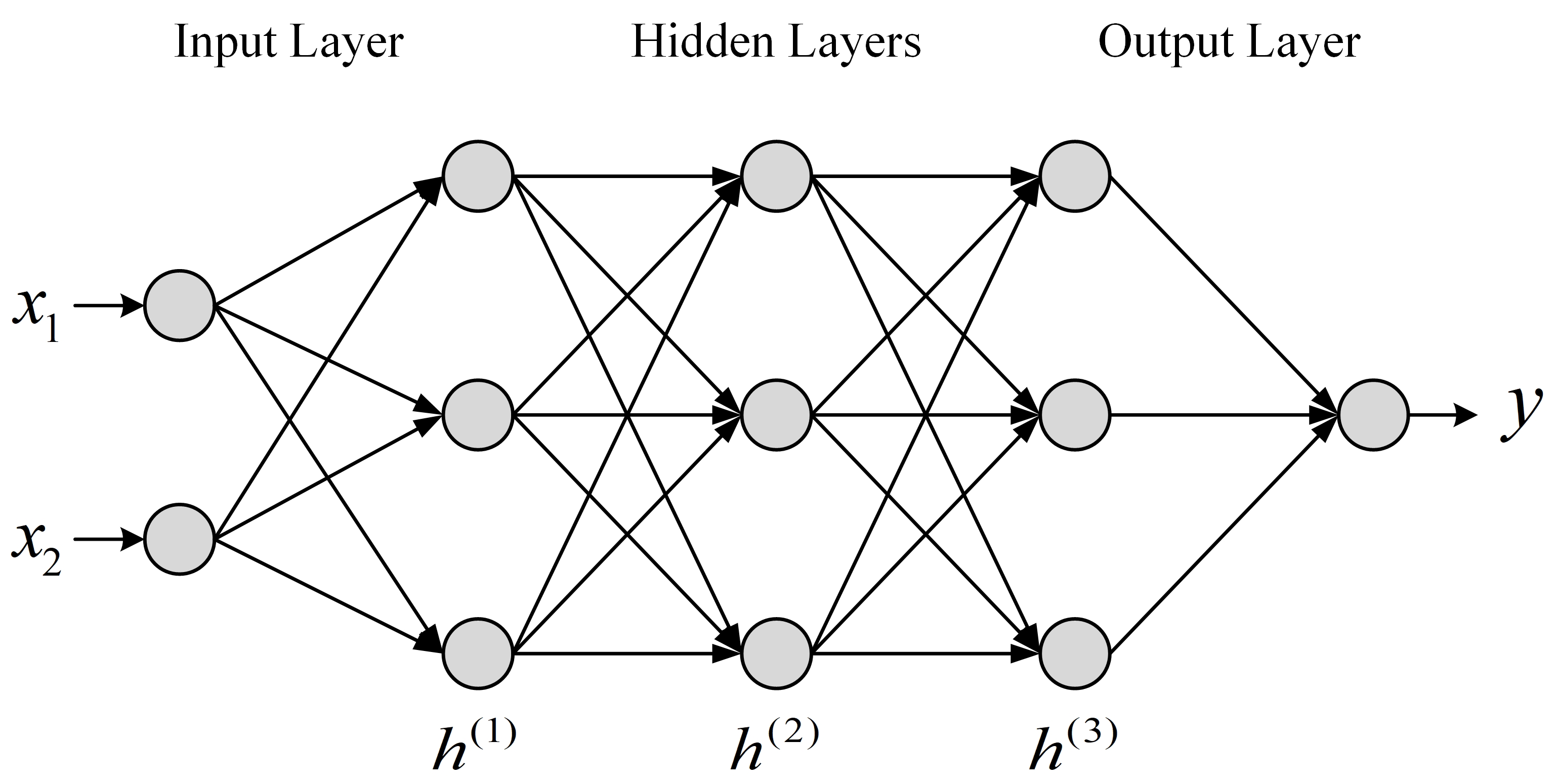}
    \caption{The framework of FNN. $x_1$, $x_2$ are the input of FNN and $y$ is the network output. $h^{(i)}$ represents hidden layer $i$.   
    }
    \label{fig2.1}
\end{figure}

A preselected neural network is used to approximate the solution. A typical network used in PINN is the FNN,  which comprises an input layer,  several hidden layers, and an output layer; see Figure \ref{fig2.1}. The input layer encodes the spatial and temporal coordinates of sampling points. The hidden layers utilize linear layers and nonlinear activation functions, such as tanh, sigmoid, and swish \cite{swish}, 
to map the linear or nonlinear relationships between inputs and outputs. The output layer generates an approximate  solution at sampling points.  

\par
Each hidden layer includes a linear layer and the activation function. In FNN, 
\begin{equation}
    \label{equ2.1}
    h^{(i)}=\sigma (w^{(i)} h^{(i-1)}+b^{(i)}).
\end{equation}
$\sigma$ is the activation function, which is generally a nonlinear function. The network parameters contain weight $w^{(i)}$ and bias $b^{(i)}$.
The FNN models the nonlinear relationship between input and output variables through multiple hidden layers. The number of layers (network depth) and the number of neurons in each layer (network width) influence its feature representation and learning capabilities.

\par
The optimization of network parameters in PINN uses the gradient descent method. The gradient descent method is an optimization algorithm that adjusts parameters in the direction opposite to the gradient of the loss function with respect to the network parameters in order to minimize the loss function and find the optimal parameters. The optimizer of PINN is typically chosen as Adam, an adaptive gradient descent optimization method.

\subsection{Search Methods}
\indent \indent 
The network hyperparameters of FNN in PINN include network width, network depth, and activation function.  The activation function is typically pre-selected,  so researchers focus primarily on determining the optimal width and depth of the network. Traditional search methods include grid search,  random search,  and Bayesian optimization. 
\par
Grid search seeks all network architectures in the search space. It is guaranteed that the optimal architecture can be found at the cost of searching over all possible networks. 
The random search method randomly samples over the search space and finds the best architecture over these samples. It has a lower computational cost but may generate suboptimal results. Bayesian optimization,  such as HyperOpt \cite{hyperopt}, 
employs a probabilistic model to predict the next hyperparameter combination to evaluate. In principle, this method is capable of finding the optimal architecture with fewer iterations than the grid search. Practically, it may perform poorly during the early searches and thus affect the final result. 
\par
Compared to traditional methods, neural architecture search (NAS) is an emerging search approach. NAS is a component of automated machine learning \cite{3-6} that aims to automatically design the optimal neural network architecture through various search strategies. It has been widely used to hyperparameter search for convolutional neural networks (CNN).
Among various methods in NAS,  gradient-based optimization methods achieve the best compromise between search time and architectural accuracy. One such method is the differentiable neural architecture search (DARTS) \cite{DARTS}.

\begin{figure}[htbp]
    \centering
    \includegraphics[width=0.3\linewidth]{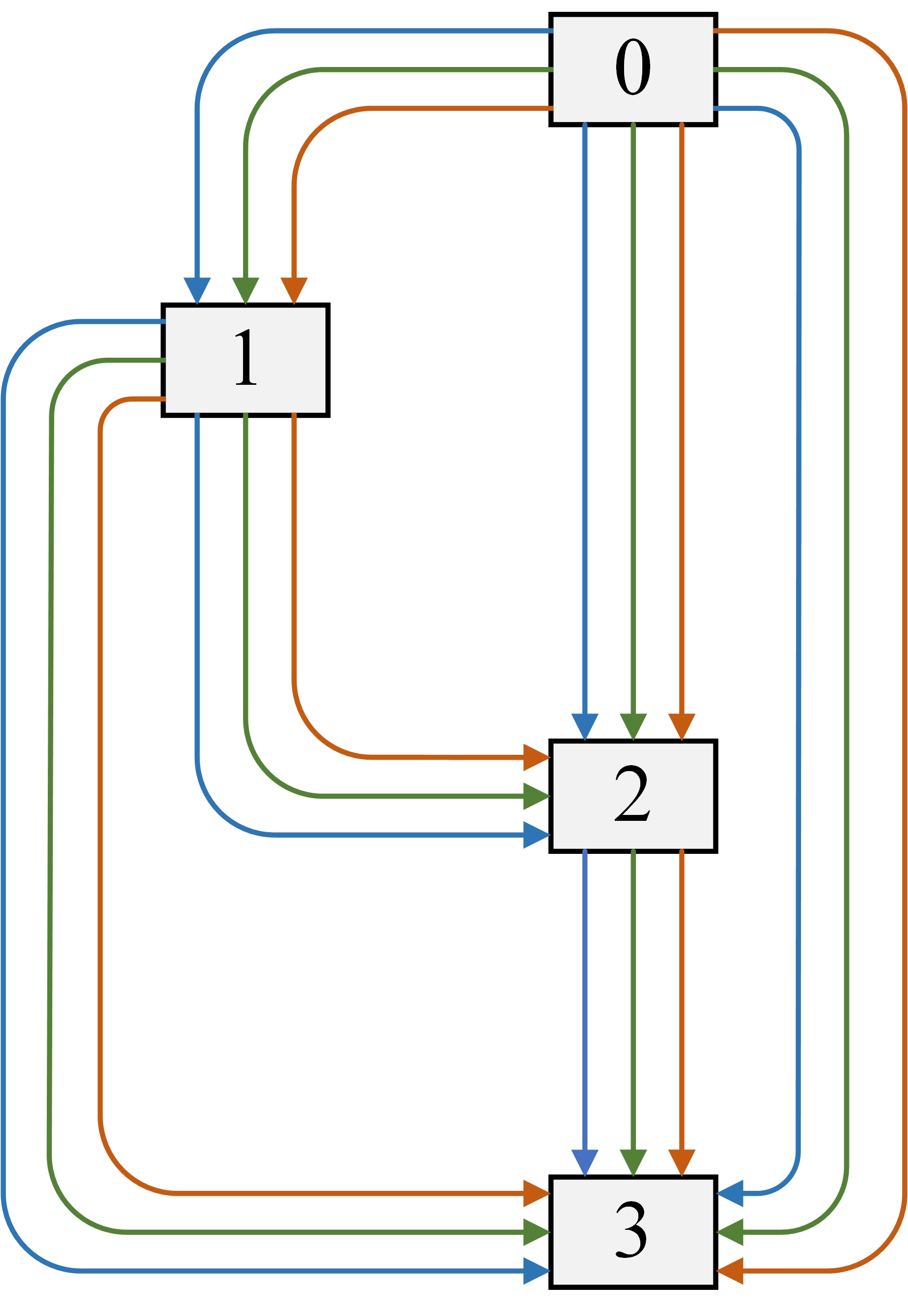}
    \caption{A cell of DARTS. All edges point from preceding nodes to subsequent nodes. There are multiple edges between two nodes, each corresponding to a candidate operation.
    }
    \label{fig2.2}
\end{figure}

DARTS seeks the structure of a cell,  which can be stacked to form the entire architecture; see Figure \ref{fig2.2}. Each cell represents as a directed acyclic graph (DAG) with $N$ nodes $\{x_i\}_{i=0}^{N-1}$, where $x_i$ is the latent representation of node $i$. Edge $(i,j)$ in the DAG represents the transformation of information from node $x_i$ to $x_j$,  containing candidate operations between $x_i$ and  $x_j$. For example,  the candidate operations of CNN include zero, skip-connection, convolution, and max-pool. All candidate operations form the operation space $\mathcal{O}$. The weight of the candidate operation is the smooth architecture parameter $\alpha^{(i,j)}$ after applying the softmax function. 
Each intermediate node $j$ is computed based on all preceding nodes: 
\begin{equation}
    \label{2.6}
    \begin{aligned}
         x_j &= \sum_{i<j} \bar{o}^{(i,j)}(x_i),  \\
        \bar{o}^{(i,j)}(x_i) &= \sum_{o \in \mathcal{O}}
        \frac{\exp(\alpha_o^{(i,j)})}{\sum_{o' \in \mathcal{O}} \exp(\alpha_{o'}^{(i,j)})}o(x) . 
    \end{aligned}
\end{equation}

\par

The output node $x_{N-1}$ is the joint representation of all intermediate nodes. For convenience, we denote the weights and biases of the DARTS network as $w$.  
The application of the softmax function transforms the original discrete search space into a continuous space, allowing the optimization of parameters using gradient-based methods. 
By updating the architecture parameters and network parameters simultaneously, DARTS achieves a shorter search time compared to other NAS methods.
We denote the loss function of the training set and validation set as ${L}_{train}$ and ${L}_{val}$, respectively. 
\par
Searching network parameters $w$ and architecture parameters $\alpha$ can be formulated as the bi-level optimization problem:
\begin{equation}
    \label{2.7}
    \begin{aligned}
         \min_{\alpha} &\quad {L}_{val}(w^{*}(\alpha),\alpha),  \\
         \text{s.t.} &\quad w^{*}(\alpha) = \arg \min_{w} {L}_{train}(w,\alpha). 
    \end{aligned}
\end{equation}
Once $\alpha$ is updated, the candidate operation $o$ with the highest weight $\alpha_o^{(i,j)}$ is selected to replace the original mixed operation for each edge $(i,j)$. After the training, the final architecture is the optimal architecture. 
\par
DARTS is typically applied to search the optimal architecture of CNN. Due to its excellent search efficiency, we apply DARTS to the optimal architecture search of PINN, aiming to obtain a new approach with better performance than traditional search methods.

\section{PINN-DARTS}
\label{section3}
\subsection{Framework}
\indent \indent
The goal of PINN-DARTS is to find the optimal network architecture of FNN in PINN for solving PDEs, specifically the number of layers and the number of neurons in each layer. We replace the hidden layers of FNN with the PINN-DARTS network and keep the input layer and output layer unchanged.\par
\begin{figure}[htbp]
    \centering
    \includegraphics[width=0.75\linewidth]{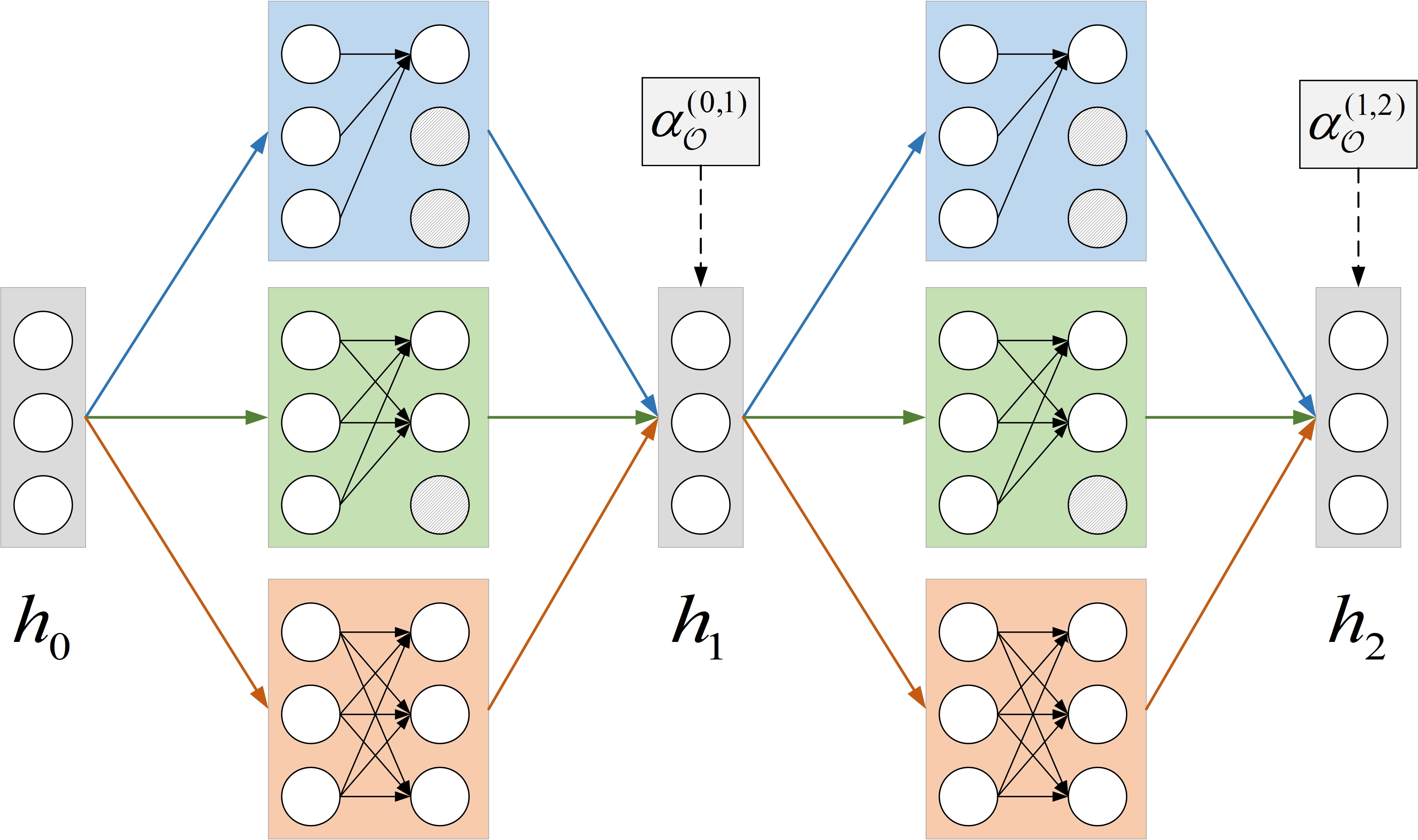}
    \caption{The structure of PINN-DARTS. Each node only has edges with its adjacent nodes. The candidate operations consist of  linear layers with different widths and the skip-connection operation. The candidate operations are weighted by $\alpha_{\mathcal{O}}^{(i-1,i)}$. $h_i$ is the output tensor of the hidden layer corresponding to node $i$, and it is the sum of weighted candidate operation inputs.  
    }
    \label{fig3.1}
\end{figure}

The structure of PINN-DARTS is shown in Figure \ref{fig3.1}. In contrast to CNN, which is typically formed by stacking multiple network blocks, FNN is composed of multiple linear layers. Therefore, the optimal architecture of a single cell in PINN-DARTS is regarded as the final network architecture.
The connection relationships between different layers in FNN are typically simple. Each network layer only has connections with its adjacent layers,  forming a linear structure rather than a graph structure.  In PINN-DARTS, each intermediate node $i$ only has edges with its parent node $i-1$ and its child node $i+1$.
\par 
The search space of PINN-DARTS is FNN with various network depths and layer widths. Differing from traditional search methods, the architectures in our consideration allow different numbers of neurons for layers in an FNN, yielding a much larger search space. The transformation between two adjacent nodes represents a linear layer. The maximum network depth of search space equals the number of PINN-DARTS nodes minus one. The candidate operations of edges contain linear layers with various widths and the skip-connection operation. The output tensors generated by linear layers with varying widths differ in shape, making direct computation impossible. To overcome this, we adopt zero-padding to standardize the tensor shapes to that of the linear layer with maximum width in search space.  

\par

PINN-DARTS network consists of $N$ nodes $\{h_i\}_{i=0}^{N-1}$.  The $h_i$ corresponding to each node $i$ represents a latent output tensor of the hidden layer.  The shape of $h_i$ is the same as the output tensor of the linear layer with maximum width. Edges exist only between adjacent nodes, represented as  $\{(i-1, i)\}_{i=1}^{N-1}$. Each edge corresponds to a candidate operation. Each candidate operation has an architecture parameter $\alpha_o^{(i-1, i)}$ as its weight.
Consistent with DARTS, architecture parameters need to be smoothed using softmax function to make the search space continuous. Each node $i$ gathers all its input edges from $i-1$ to form its representation $h_i$.
\begin{equation}
    \label{3.1}
    \begin{aligned}
    h_i &= \sum_{o \in \mathcal{O}} \text{softmax}(\alpha_o^{(i-1,i)}) o(h_{i-1}) \\      
    &= \sum_{o \in \mathcal{O}} 
    \frac{\exp(\alpha_o^{(i-1,i)})}{\sum_{o' \in \mathcal{O}} \exp(\alpha_{o'}^{(i-1,i)})} o(h_{i-1}).  \\
    \end{aligned}
\end{equation}
\par
PINN-DARTS includes a searching phase and an evaluation phase. In the searching phase, architecture parameters $\alpha$ and network parameters $w$ are updated simultaneously using the gradient descent method. The detailed update formula is provided in Section \ref{subsection3.2}. After the searching phase is completed, the candidate operation with the largest $\alpha_o^{(i-1, i)}$ is selected as the optimal operation between nodes, resulting in the discretized optimal network architecture. In the evaluation phase, the network architecture is fixed to be the optimal architecture obtained during the searching phase. The network parameters are trained to obtain the numerical solution of PDE and calculate the error with respect to the exact solution. 
\par
It is found that DARTS tends to contain excessive skip-connection operations in the optimal network architecture. At the continuous level, the skip-connection operation forms a structure similar to ResNet, showing an advantage in competition with other candidate operations. However, after switching to the discrete setup,  the inclusion of skip-connection operation results in the deterioration of network performance. We also observed the phenomenon in our early experiments for PINN-DARTS.
\par
To solve this problem, we incorporate DARTS+ \cite{DARTS+} and FairDARTS \cite{FairDARTS} into our work. DARTS+ adds an early stopping mechanism to prevent excessive training. PINN-DARTS+ refers to PINN-DARTS with the early stopping mechanism. In PINN-DARTS+, the early stopping condition is set as the $\alpha$ ranking remains unchanged in 100 consecutive iterations or the optimal network architecture contains 6 or more skip-connection operations. 
FairDARTS employs the sigmoid function in place of softmax function to prevent competition among candidate operations. PINN-DARTS with the sigmoid function is named as PINN-SDARTS. The variant that incorporates both improvements is called PINN-SDARTS+. 
\par
Algorithm \ref{algo1} summarizes the algorithmic details of PINN-DARTS.
\begin{algorithm}[htbp]
    \caption{PINN-DARTS}
    \label{algo1}
    \begin{algorithmic}
        \STATE Distribute candidate operation ${o}^{(i-1,i)}$ parametrized by $\alpha_o^{(i-1,i)}$ for each edge $(i-1,i)$;
        \STATE Node $i$ gathers all input edges $(i-1,i)$ and 
        calculate $h_i$;
        \IF{SDARTS}
            \STATE $h_i=\sum \text{sigmoid}(\alpha^{(i-1,i)}) o(h_{i-1})$;
        \ELSE
            \STATE $h_i=\sum \text{softmax}(\alpha^{(i-1,i)}) o(h_{i-1})$;
        \ENDIF
        \WHILE{not converge}
            \STATE Update architecture parameters and network parameters:
            \STATE $\alpha = \alpha - \nabla_{\alpha} {L}_{test}(w,\alpha)$;
            \STATE $w = w - \nabla_{w} {L}_{train}(w,\alpha)$;
            \IF{DARTS+ and satisfy condition}
                \STATE Break loop and early stopping;
            \ENDIF
        \ENDWHILE
        \STATE Select optimal operation based on $\alpha$ and output the optimal network architecture.
    \end{algorithmic}
\end{algorithm}

\subsection{Update Details}
\label{subsection3.2}
\indent \indent
For real applications,  the exact solution to a given PDE is typically not available,  thus, a method that does not require any prior information of the exact solution is favorable.
We point out that PINN-DARTS is fully unsupervised over the entire searching phase. An update of the architecture parameters $\alpha$ only requires the information from the test set, and an update of the network parameters $w$ requires information from the training set. Precisely, the update of the network parameters $w$ uses the PINN loss of the training set, and the update of the architecture parameters $\alpha$ uses the PINN loss of the test set.  $M$ collocation points for the training set and $T$ collocation points for the test set are employed, respectively. ${L}_{train}$ and ${L}_{test}$ are the PINN losses corresponding to the training set and the test set.
\begin{equation} 
    \label{3.2.1}
    \begin{aligned}
        {L}_{train} & =\frac{1}{M_f} \sum_{i=1}^{M_f} \left| \mathcal{L}(u(\boldsymbol{x}_i,t_i)) \right|^2+
        \frac{1}{M_b} \sum_{i=1}^{M_b} \left| u(\boldsymbol{x}_i,t_i) - b(\boldsymbol{x}_i,t_i) \right|^2 +\frac{1}{M_i} \sum_{i=1}^{M_i} \left| u(\boldsymbol{x}_i,0) - u_0(\boldsymbol{x}_i) \right|^2 , \\
        {L}_{test} & =\frac{1}{T_f} \sum_{i=1}^{T_f} \left| \mathcal{L}(u(\boldsymbol{x}_i,t_i)) \right|^2+
        \frac{1}{T_b} \sum_{i=1}^{T_b} \left| u(\boldsymbol{x}_i,t_i) - b(\boldsymbol{x}_i,t_i) \right|^2 +\frac{1}{T_i} \sum_{i=1}^{T_i} \left| u(\boldsymbol{x}_i,0) - u_0(\boldsymbol{x}_i) \right|^2 . \\
    \end{aligned}
\end{equation}
The architecture parameters $\alpha$ and network parameters $w$ are updated using the gradient descent method:
\begin{equation}
    \label{3.2}
    \begin{aligned}
        \alpha &=\alpha- \nabla_{\alpha} {L}_{test}(w,\alpha),  \\ 
        w &=w- \nabla_{w} {L}_{train}(w,\alpha) . 
    \end{aligned}
\end{equation}
\par
The success of PINN-DARTS relies heavily on its unsupervised nature. In what follows, we shall demonstrate that the PINN loss and the relative $L^2$ error between the numerical solution associated with a given architecture and the exact solution
are strongly correlated. Consider the Poisson equation with Dirichlet boundary conditions in \eqref{3.3}. The exact solution is set to be $u(x,y) = \cos(\pi x)\cos(\pi y)$ (simple) and $u(x,y) = \cos(2\pi x)\cos(2\pi y)$ (complex). 
\begin{equation}
    \label{3.3}
    \left\{ 
        \begin{aligned}
         -\Delta u(x,y) &= f(x,y) \quad & x, y \in \Omega,  \\ 
         u(x,y) &= g(x,y) \quad & x, y \in \partial \Omega. 
        \end{aligned}
    \right. 
\end{equation}

\begin{figure}[htbp]
    \centering
    \subfigure[Simple]{\includegraphics[width=0.49\linewidth]{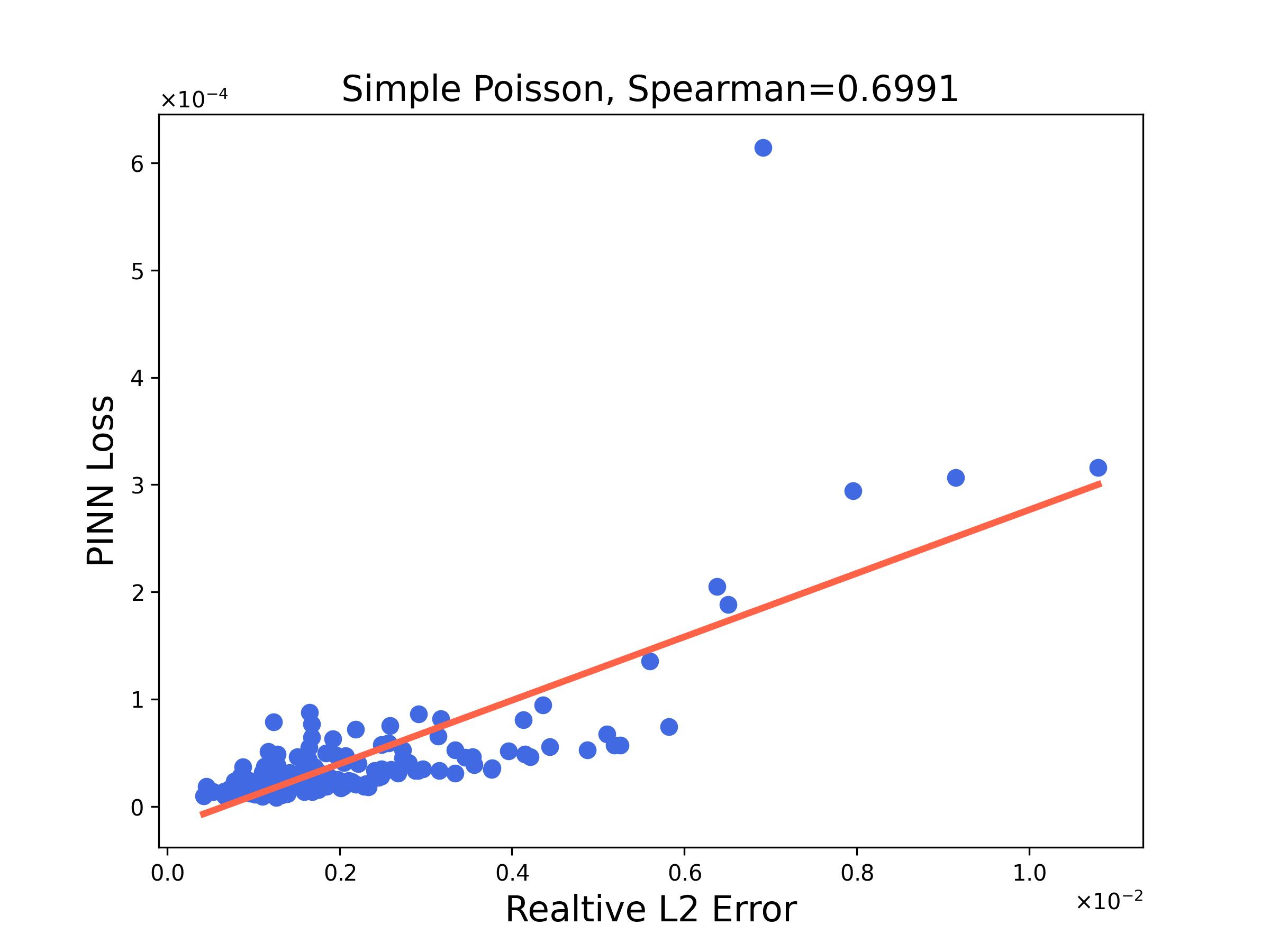}}
    \subfigure[Complex]{\includegraphics[width=0.49\linewidth]{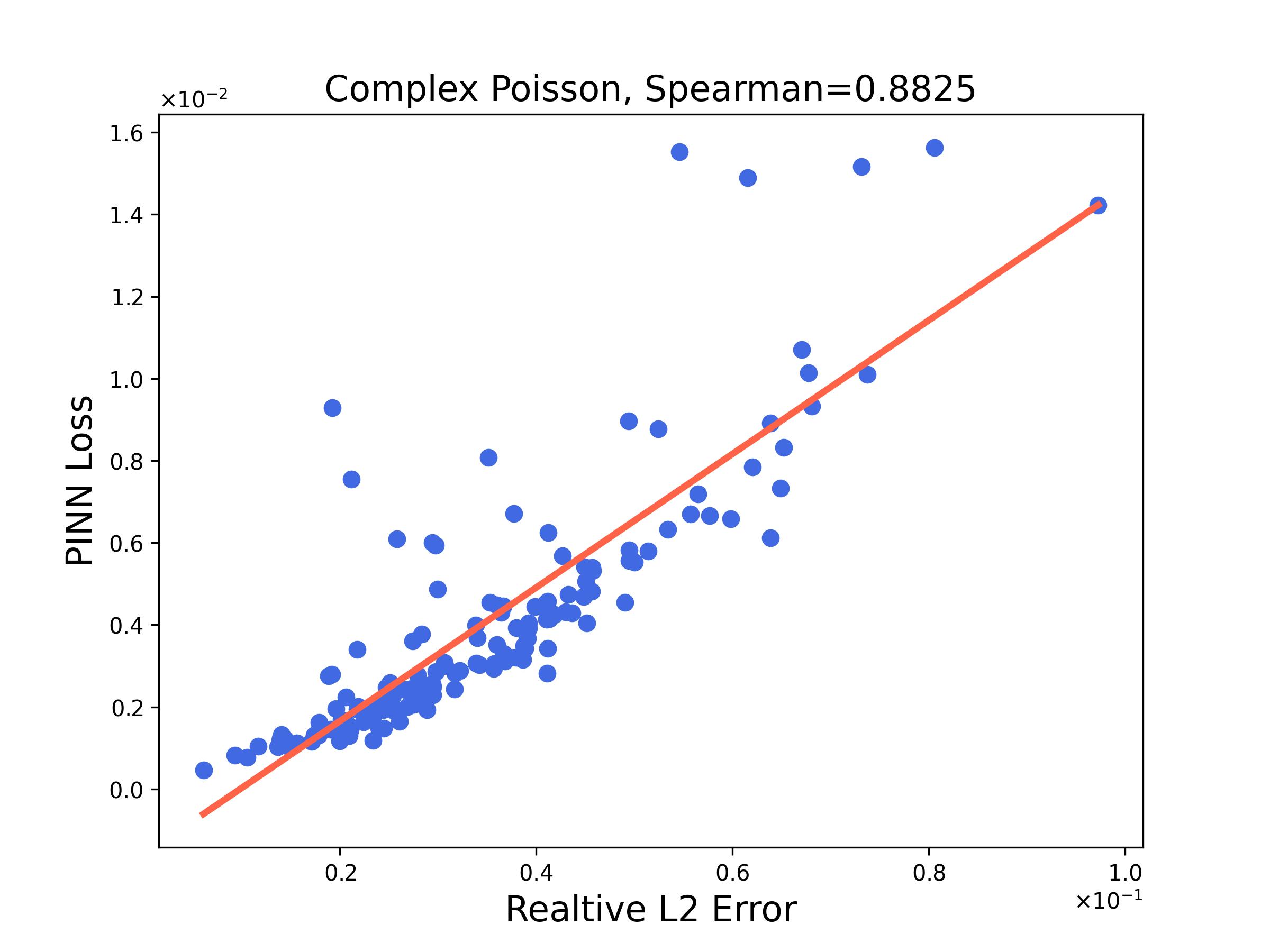}}
    \caption{Spearman's rank correlation coefficients of PINN loss and relative $L^2$ error for a series of architectures with different widths and depths for the Poisson equation with simple and complex solutions.}
    \label{fig3.2}
\end{figure}

To facilitate the discussion, we only consider FNN with network widths $\{100, 200, 300, 400\}$ and network depths between 1 and 8. 
The total 32 architectures are trained with 5 random seeds, resulting in 160 pairs of data (PINN loss, relative $L^2$ error). To analyze this data set, we apply Spearman's rank correlation coefficient, which converts the values of data points into ranks and evaluates the correlation between two variables by comparing their ranks. Scatter plots and correlation coefficients with simple and complex exact solutions are shown in Figure \ref{fig3.2}. The corresponding Spearman's rank correlation coefficients are 0.6991 and 0.8825, indicating the strong correlation between PINN loss and the relative $L^2$ error for a given architecture.

\section{Numerical Experiments}
\label{section4}
\indent \indent
In this section,  we present numerical experiments to demonstrate the superiority of PINN-DARTS methods.  Examples include the Poisson equation, heat equation, wave equation, and Burgers' equation.
Specifically, we compare PINN-DARTS methods with traditional methods. PINN-DARTS methods include PINN-DARTS,  PINN-DARTS+,  PINN-SDARTS, and PINN-SDARTS+. Traditional search methods include grid search, random search, and Bayesian search. 
We choose the TPE search of HyperOpt \cite{hyperopt} for the Bayesian search. Grid search explores all architectures (with even widths) in the search space and thus outputs the global optimal network architecture. 
It is found that the optimal network architecture relies on both the PDE and the solution complexity. Detailed optimal architectures are provided in Appendix \ref{app1}.
\par
The number of iterations in the random and Bayesian search is set to be 2 or 5, ensuring that their running times are close to those of the PINN-DARTS methods. We compare these methods in terms of search time and architectural accuracy. To enhance the reliability of the experimental results, we repeat our experiments under 5 different network parameter initializations. The training set and the test set remain fixed across different initializations. The search results of 5 initializations are labeled as Init1 to Init5. We change the network parameter initialization by varying the random seed. 
The search space is defined as FNN with up to 8 layers, and each layer contains 100, 200, 300, or 400 neurons. 
\par
We evaluate the performance of a method using search time and error ratio. The search time is defined as the running time of the searching phase and evaluation phase measured in minutes. The error ratio is defined as the average architecture error of the current search method divided by the average error of the grid search. All experiments are conducted on 1 NVIDIA A40 GPU.

\subsection{Poisson Equation}
\indent \indent
Consider a Poisson equation with Dirichlet boundary condition in \eqref{4.1} over $\Omega = [0,1] \times [0,1]$. To evaluate the performance of PINN-DARTS with different exact solutions,  we select two exact solutions: a simple solution $u(x,y) = \cos(\pi x)\cos(\pi y)$ and a complex solution $u(x,y) = \cos(2\pi x)\cos(2\pi y)$. 
\begin{equation}
    \label{4.1}
    \left\{ 
        \begin{aligned}
         -\Delta u(x,y) &= g(x,y) \quad &x,y \in \Omega,  \\
         u(x,y) &= f(x,y) \quad &x,y \in \partial \Omega. 
        \end{aligned}
    \right. 
\end{equation}
\par
The training set consists of 5,000 random points sampled inside the domain and 200 random points from the boundary. The test set comprises uniformly sampled points of the $100 \times 100$ grid. For traditional search methods and the evaluation phase of PINN-DARTS,  the initial learning rate of network parameters is $5 \times 10^{-4}$ and decays to $1 \times 10^{-5}$ as the training process proceeds. The total number of iterations is set to 10,000. The searching phase for PINN-DARTS is limited to a maximum of 2,000 iterations, and early stopping is applied if the convergence condition is satisfied. During the searching phase, the learning rate of  architecture parameters is set to be $2 \times 10^{-1}$, and the learning rate of network parameters is set to be $1 \times 10^{-4}$. For both phases, Adam is employed as the optimizer, and tanh is used as the activation function. The same setup in PINN-DARTS is used for other PINN-DARTS variants. In Table \ref{tab4.1} and Table \ref{tab4.2}, we record the average search time and the error ratio of different search methods. 

\begin{table}[htbp]
    \setlength{\belowcaptionskip}{0.1cm}
    \caption{Optimal architecture results of Poisson equation with the simple solution. Relative $L^2$ error, search time, and error ratio are averaged over results with 5 initializations.}
    \label{tab4.1}
    \centering
    \resizebox{0.7\textwidth}{!}{
    \begin{tabular}{cccc}
        \toprule
        Search Method & Relative $L^2$ Error & Search Time (min)  & Error Ratio \\
        \midrule
        Grid Search  & 6.532e-04 &206 & 100.0\%\\
        Random Search, Iter=5  & 8.877e-04 &34 & 135.9\%  \\
        Random Search, Iter=2  & 1.482e-03 &12 & 226.9\%  \\
        Bayesian Search, Iter=5   & 9.429e-04 &31 & 144.3\% \\
        Bayesian Search, Iter=2   & 1.425e-03 &13 & 218.2\% \\
        PINN-DARTS   & 1.389e-03 & 30 & 212.6\% \\
        PINN-DARTS+   & 1.389e-03  & 10 & 212.6\% \\
        PINN-SDARTS   & 1.260e-03 & 28 & 192.9\% \\
        PINN-SDARTS+  & 1.308e-03 & 12 & 200.3\% \\
        \bottomrule
    \end{tabular}}
\end{table}

\begin{table}[htbp]
    \setlength{\belowcaptionskip}{0.1cm}
    \caption{Optimal architecture results of Poisson equation with the complex solution.
    Relative $L^2$ error, search time, and error ratio are averaged over results with 5 initializations.
    }
    \label{tab4.2}
    \centering
    \resizebox{0.7\textwidth}{!}{
    \begin{tabular}{cccc}
        \toprule
        Search Method & Relative $L^2$ Error & Search Time (min) & Error Ratio\\
        \midrule
        Grid Search  & 1.150e-02 & 316 & 100.0\% \\
        Random Search, Iter=5  & 2.203e-02 & 55 & 191.6\%  \\
        Random Search, Iter=2  & 2.595e-02 & 22 & 225.7\% \\
        Bayesian Search, Iter=5  & 1.712e-02 & 48 & 148.9\% \\
        Bayesian Search, Iter=2  & 1.896e-02 & 20 & 164.9\% \\
        PINN-DARTS   & 1.784e-02 & 35 & 155.1\% \\
        PINN-DARTS+  & 1.522e-02 & 15 & 132.4\%   \\
        PINN-SDARTS   & 4.024e-02 & 35 & 350.0\% \\
        PINN-SDARTS+   & 2.867e-02 & 14 & 249.4\% \\
        \midrule
    \end{tabular}}
\end{table}

\par
\begin{figure}[htbp]
    \centering
    \subfigure[Simple]{
    \label{sub4.1.1}
    \includegraphics[width=0.45\linewidth]{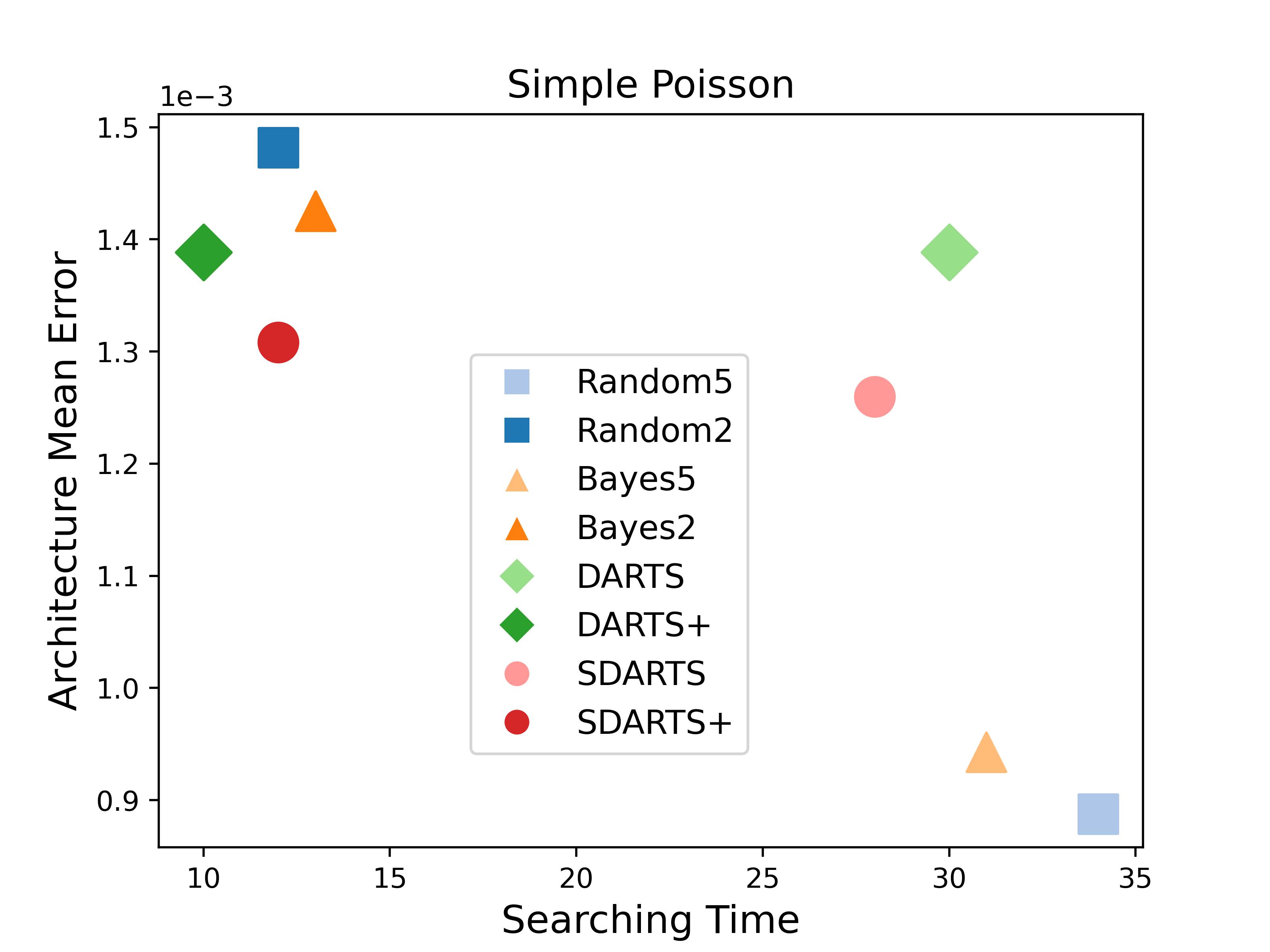}}
    \subfigure[Complex]{
    \label{sub4.1.2}
    \includegraphics[width=0.45\linewidth]{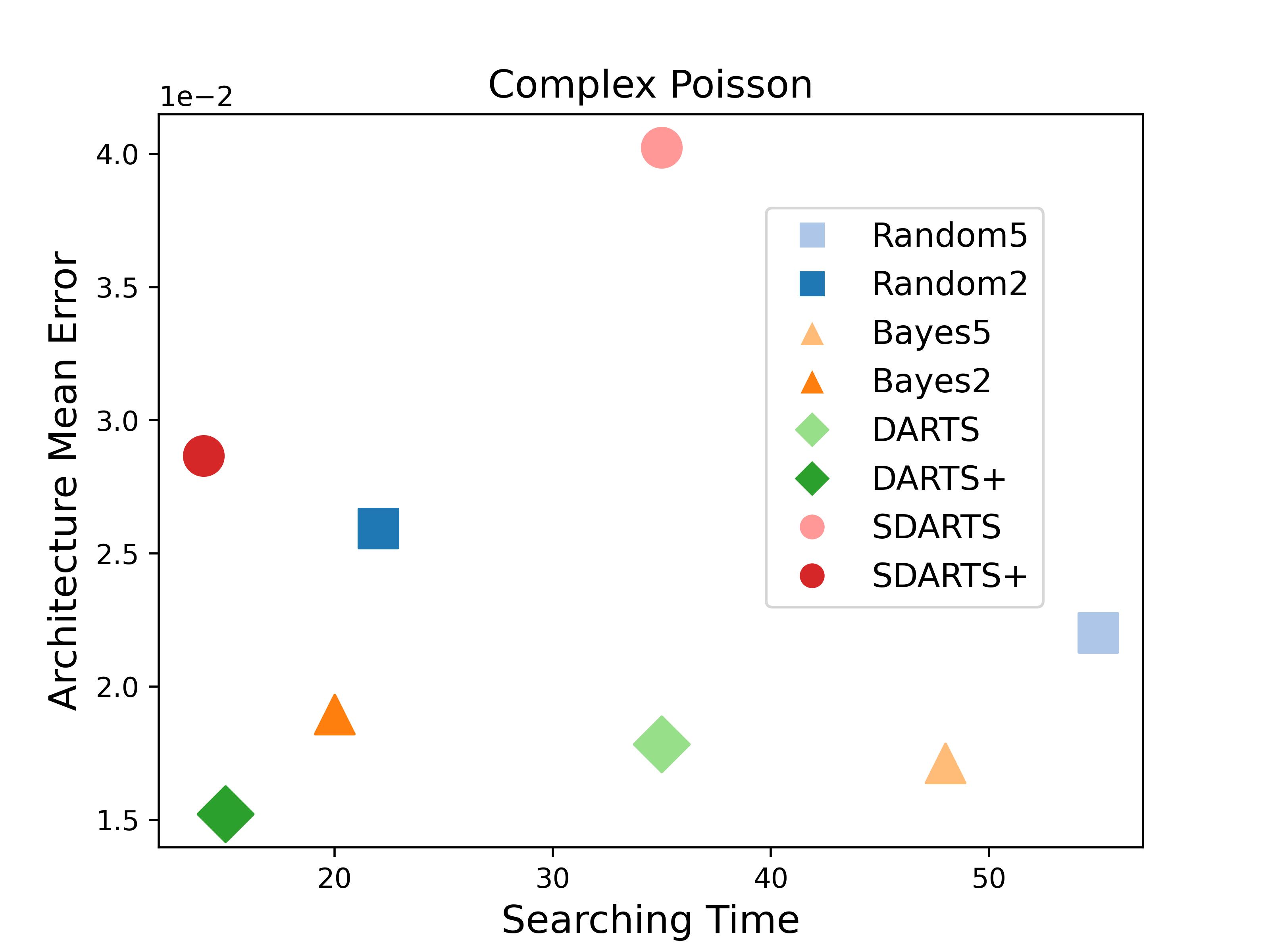}}
    \caption{Comparison of PINN-DARTS methods and traditional methods for the Poisson equation with simple and complex solutions.}
    \label{fig4.1}
\end{figure}

For the simple solution case, PINN-SDARTS+ achieves the best result among these methods within a short search time. However, if there is no limit on the search time, traditional methods win; see Figure \ref{sub4.1.1}. The performance of PINN-DARTS and PINN-DARTS+ improves significantly for the complex solution case. The architecture error of PINN-DARTS+ outperforms other methods, while its search time remains the shortest; see Figure \ref{sub4.1.2}. As the solution complexity increases, PINN-DARTS and PINN-DARTS+ perform better,  
indicating the advantage of architectures with uneven widths, possibly due to their stronger feature representation and learning capabilities. 

\subsection{Heat Equation}
\indent \indent
Consider the heat equation in \eqref{4.2} with simple and complex exact solutions: $u(x,t) = e^{-t} \sin(\pi x)$ and $u(x,t) = e^{-t} \sin(2\pi x)$.
\begin{equation}
    \label{4.2}
    \left\{ 
        \begin{aligned}
         \frac{\partial u(x,t)}{\partial t} - \frac{\partial^2 u(x,t)}{\partial x^2} &= g(x,t) \quad &x \in [-1,1],\ t \in [0,1],  \\
         u(-1,t) = u(1,t) &= f(t),  \\
         u(x,0) &= u_0(x). 
        \end{aligned}
    \right. 
\end{equation}
\par
The training set includes 5,000 randomly sampled points inside the domain and 200 randomly sampled points from the boundary.  
The test set includes $200 \times 100$ uniformly sampled points.  
For the searching phase,  the learning rate of architecture parameters and network parameters are set to $2 \times 10^{-1}$ and $1 \times 10^{-4}$, respectively.  
Other parameters remain the same as those used in the Poisson equation. 
All results are recorded in Table \ref{tab4.3} and Table \ref{tab4.4}.

\begin{table}[htbp]
    \setlength{\belowcaptionskip}{0.1cm}
    \caption{Optimal architecture results of heat equation with the simple solution. Relative $L^2$ error, search time, and error ratio are averaged over results with 5 initializations.
    }
    \label{tab4.3}
    \centering
    \resizebox{0.7\textwidth}{!}{
    \begin{tabular}{cccc}
        \toprule
        Search Method & Relative $L^2$ Error & Search Time (min) & Error Ratio\\
        \midrule
        Grid Search  & 8.817e-04 &177 & 100.0\% \\
        Random Search, Iter=5  & 1.141e-03 &29 & 129.5\%  \\
        Random Search, Iter=2  & 1.512e-03 &11 & 171.4\%  \\
        Bayesian Search, Iter=5   & 1.184e-03 &26 & 134.3\% \\
        Bayesian Search, Iter=2   & 1.345e-03 &10 & 152.5\% \\
        PINN-DARTS   & 1.095e-03 & 28 & 124.2\% \\
        PINN-DARTS+    & 1.095e-03 &9 & 124.2\%  \\
        PINN-SDARTS   & 1.602e-03 & 27 & 181.7\% \\
        PINN-SDARTS+   & 1.602e-03 & 10 & 181.7\% \\
        \bottomrule
    \end{tabular}}
    
\end{table}

\begin{table}[htbp]
    \setlength{\belowcaptionskip}{0.1cm}
    \caption{Optimal architecture results of heat equation with the complex solution. Relative $L^2$ error, search time, and error ratio are averaged over results with 5 initializations.
    }
    \label{tab4.4}
    \centering
    \resizebox{0.7\textwidth}{!}{
    \begin{tabular}{cccc}
        \toprule
        Search Method & Relative $L^2$ Error & Search Time (min) & Error Ratio\\
        \midrule
        Grid Search  & 1.943e-03 & 258 & 100.0\% \\
        Random Search, Iter=5  & 2.831e-03 & 32 & 145.7\% \\
        Random Search, Iter=2  & 3.840e-03 & 14 & 197.6\%\\
        Bayesian Search, Iter=5   & 3.890e-03 & 38 & 200.2\% \\
        Bayesian Search, Iter=2  & 5.017e-03 & 13 & 258.1\% \\
        PINN-DARTS  & 3.958e-03 & 27 & 203.7\% \\
        PINN-DARTS+  & 3.754e-03 & 11 & 193.2\% \\
        PINN-SDARTS & 7.893e-03 & 23 & 406.1\% \\
        PINN-SDARTS+ & 1.286e-02 & 9 & 661.8\%\\
        \bottomrule
    \end{tabular}}
    
\end{table}

\par
\begin{figure}[htbp]
    \centering
    \subfigure[Simple]{
    \label{sub4.2.1}
    \includegraphics[width=0.45\linewidth]{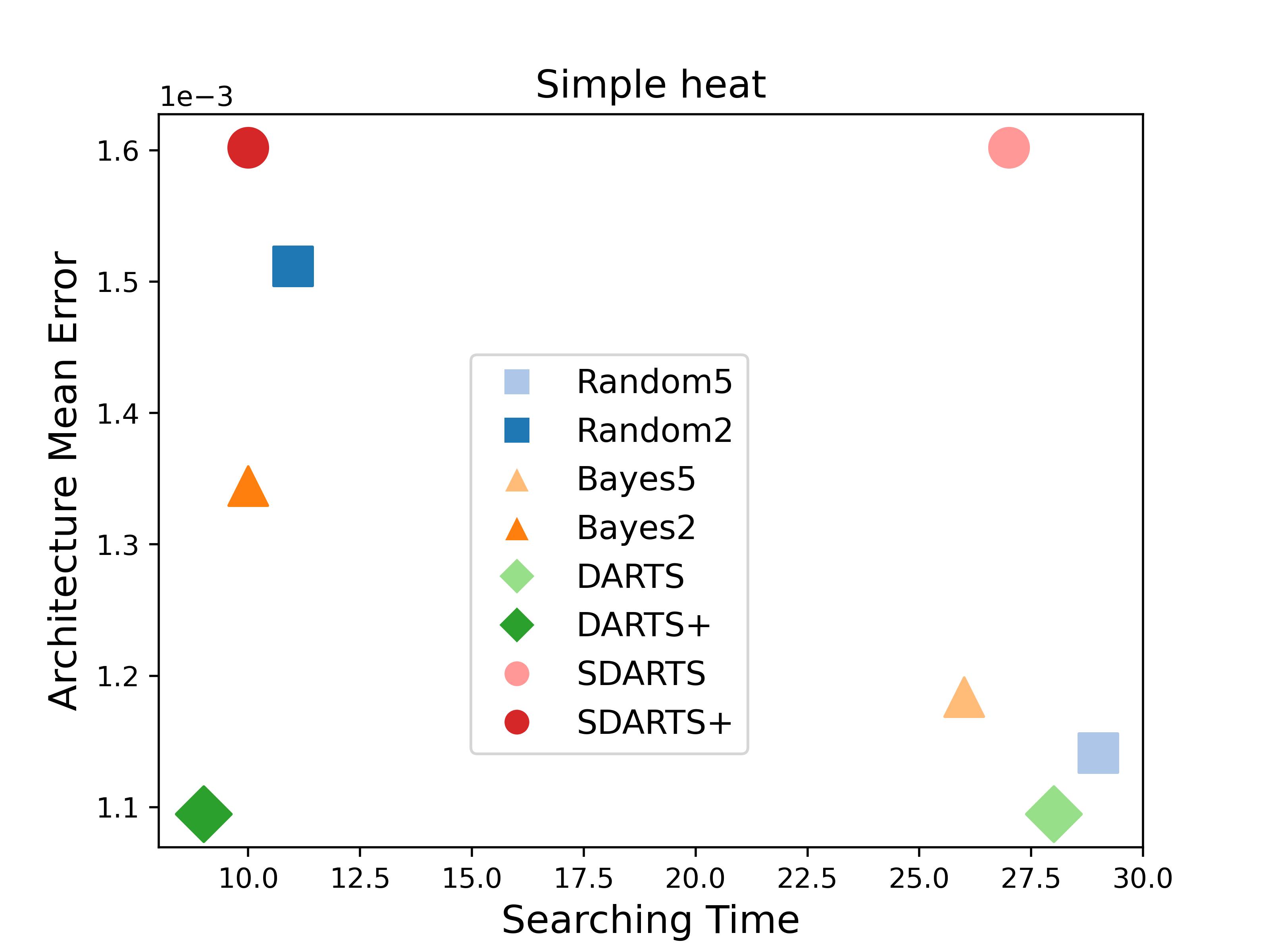}}
    \subfigure[Complex]{
    \label{sub4.2.2}
    \includegraphics[width=0.45\linewidth]{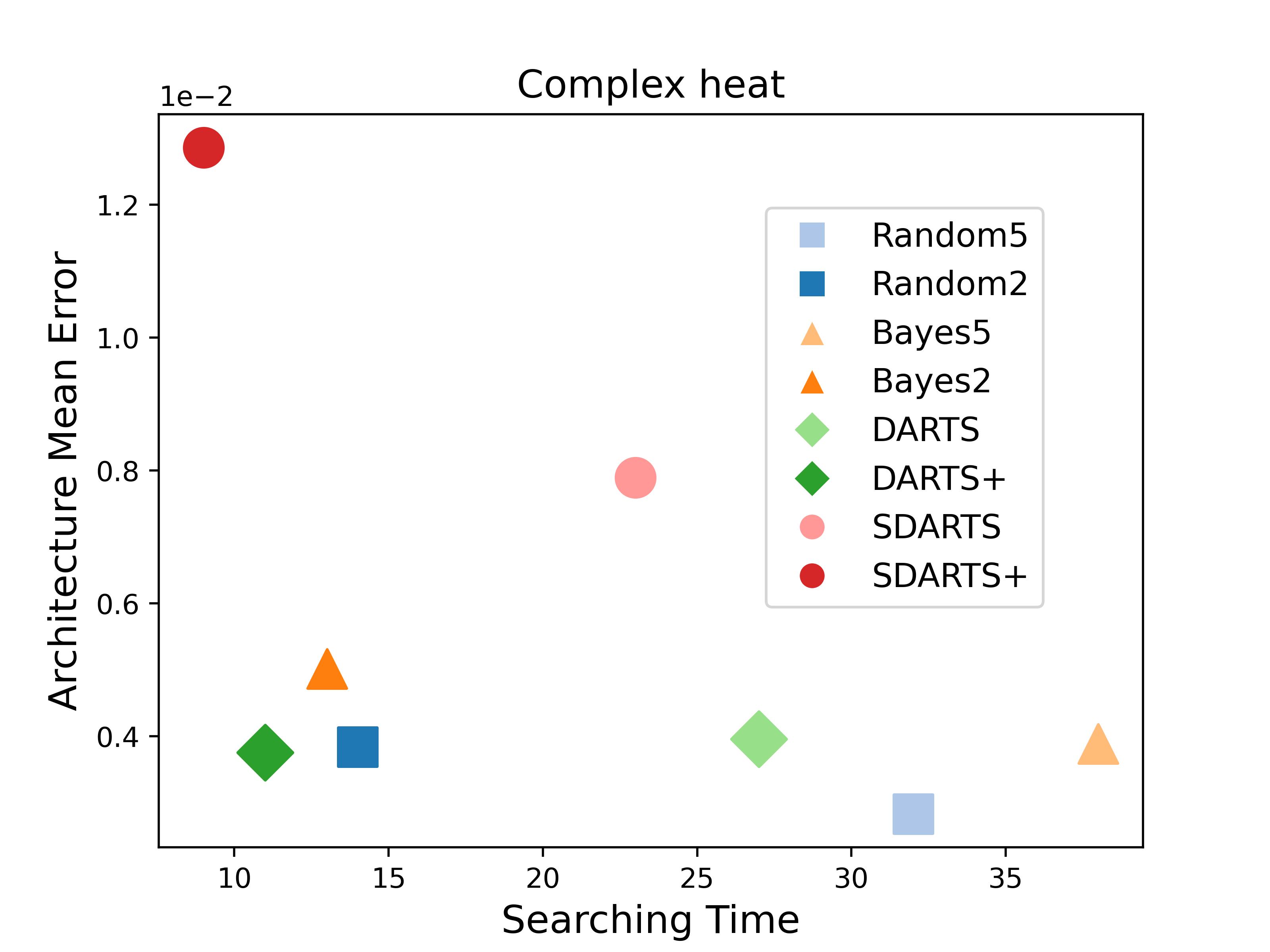}}
    \caption{Comparison of PINN-DARTS methods and traditional methods for the heat equation with simple and complex solutions.}
    \label{fig4.2}
\end{figure}

For the simple solution case, PINN-DARTS and PINN-DARTS+ outperform other search methods in Figure \ref{sub4.2.1}, while the opposite is observed for the complex solution case. From Appendix \ref{app1}, we find that the optimal network architecture (with even widths) of simple solution is (width=400, depth=6), and the optimal architecture of complex solution is (width=300, depth=2). With the previous analysis of experiment results, we conclude that PINN-DARTS and PINN-DARTS+ perform better for problems where the optimal network architecture is deeper and wider.

\subsection{Wave Equation}
\indent \indent
Consider the wave equation in \eqref{4.3} with a simple solution $u(x,t) = \sin(xt)$ and a complex solution $u(x,t) = \sin(2xt)$. 
\begin{equation}
    \label{4.3}
    \left\{ 
        \begin{aligned}
        & \frac{\partial^2 u(x,t)}{\partial t^2} - \frac{\partial^2 u(x,t)}{\partial x^2} = g(x,t) \quad x \in [0,1],  \ t \in [0,1],  \\
        & u(0,t) = f_0(t),  \ u(1,t) = f_1(t),  \\
        & u(x,0) = u_0(x),  \ \frac{\partial u(x,t)}{\partial t} \Big|_{t=0} = u_1(x). 
        \end{aligned}
    \right. 
\end{equation}
\par 
The training set contains 5,000 randomly sampled points inside the domain and 200 randomly sampled points from the boundary.  
The test set contains $100 \times 100$ uniformly sampled points.  
For the searching phase,  the learning rate of architecture parameters is set to be $5 \times 10^{-2}$,  
and the learning rate of network parameters is set to be $1 \times 10^{-4}$. Following Auto-PINN \cite{Auto-PINN}, we choose Swish as
the activation function for both phases. Other parameters remain the same as before. Table \ref{tab4.5} and  Table \ref{tab4.6} contain all optimal architecture results.

\begin{table}[htbp]
    \setlength{\belowcaptionskip}{0.1cm}
    \caption{Optimal architecture results of the wave equation with the simple solution. 
    Relative $L^2$ error, search time, and error ratio are averaged over results with 5 initializations.
    }
    \label{tab4.5}
    \centering
    \resizebox{0.7\textwidth}{!}{
    \begin{tabular}{cccc}
        \toprule
        Search Method & Relative $L^2$ Error & Search Time (min) & Error Ratio\\
        \midrule
        Grid Search  & 4.075e-04 &297 & 100.0\%\\
        Random Search, Iter=5  & 7.598e-04 &32 & 186.5\%  \\
        Random Search, Iter=2  & 1.244e-03 &12 & 305.3\%  \\             
        Bayesian Search, Iter=5   & 6.338e-04 &30 & 155.5\% \\      
        Bayesian Search, Iter=2   & 1.266e-03 &10 & 310.8\% \\
        PINN-DARTS   & 9.472e-04 & 26 & 232.5\% \\
        PINN-DARTS+    & 9.167e-04 &9 & 225.0\%  \\
        PINN-SDARTS   & 1.518e-03 & 28 & 372.4\% \\
        PINN-SDARTS+   & 2.361e-03 & 11 & 579.5\% \\
        \bottomrule
    \end{tabular}}
\end{table}

\begin{table}[htbp]
    \setlength{\belowcaptionskip}{0.1cm}
    \caption{Optimal architecture results of the wave equation with the complex solution. 
    Relative $L^2$ error, search time, and error ratio are averaged over results with 5 initializations.
    }
    \label{tab4.6}
    \centering
    \resizebox{0.7\textwidth}{!}{
    \begin{tabular}{cccc}
        \toprule
        Search Method & Relative $L^2$ Error & Search Time (min) & Error Ratio\\
        \midrule
        Grid Search  & 1.815e-03 & 334 & 100.0\%\\
        Random Search, Iter=5  & 2.512e-03 & 38 & 138.4\%  \\
        Random Search, Iter=2  & 4.693e-03 & 14 & 258.6\%  \\
        Bayesian Search, Iter=5   & 3.007e-03 & 31 & 165.7\% \\
        Bayesian Search, Iter=2   & 7.053e-03 & 13 & 388.7\% \\
        PINN-DARTS   & 4.884e-03 & 24 & 269.2\% \\
        PINN-DARTS+  & 3.444e-03 & 11 & 189.8\%   \\
        PINN-SDARTS   & 1.367e-02 & 22 & 753.6\% \\
        PINN-SDARTS+   & 1.367e-02 & 13 & 753.6\%  \\
        \bottomrule
    \end{tabular}}
\end{table}

\par
\begin{figure}[htbp]
    \centering
    \subfigure[Simple]{
    \label{sub4.3.1}
    \includegraphics[width=0.45\linewidth]{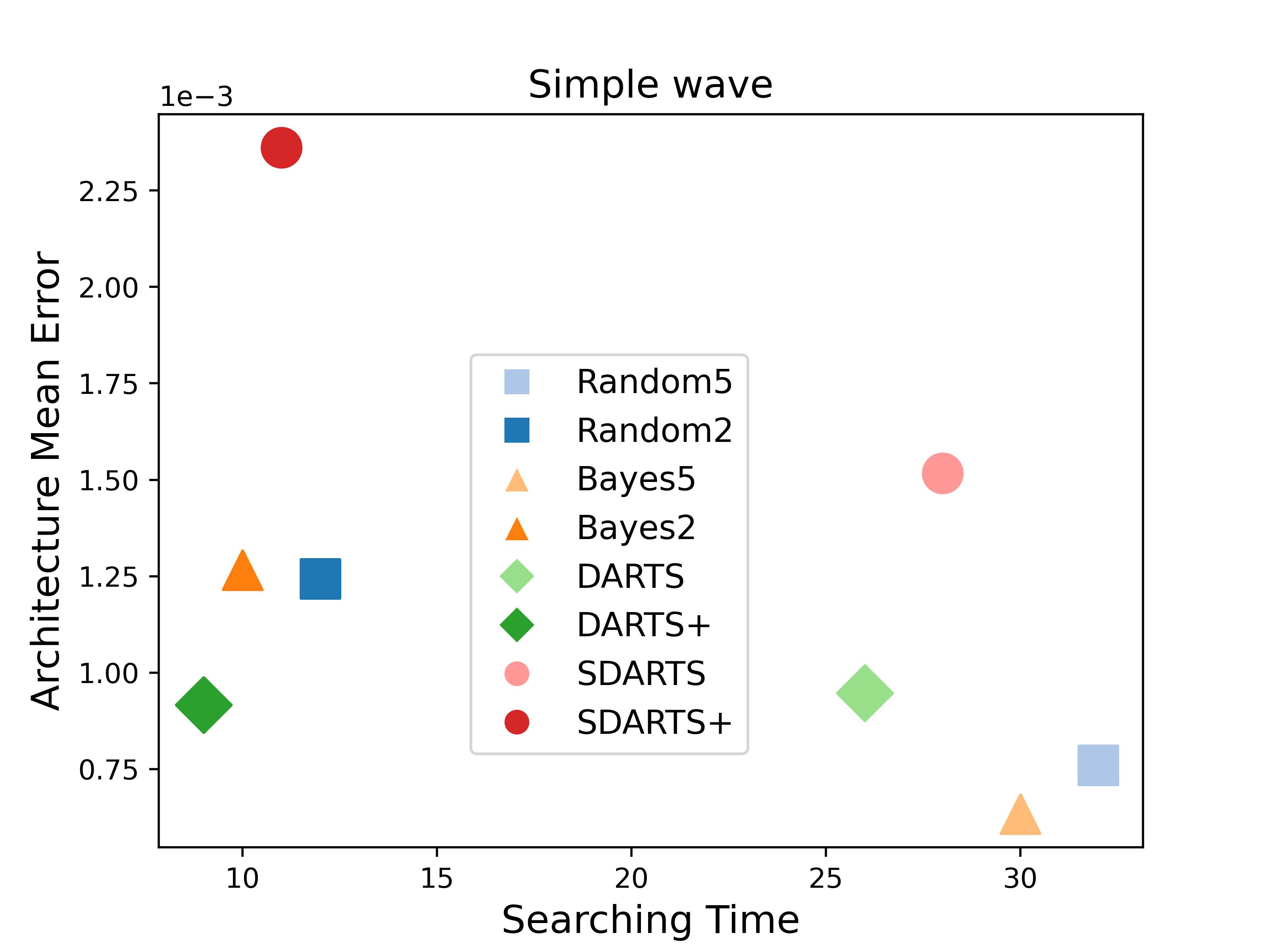}}
    \subfigure[Complex]{
    \label{sub4.3.2}
    \includegraphics[width=0.45\linewidth]{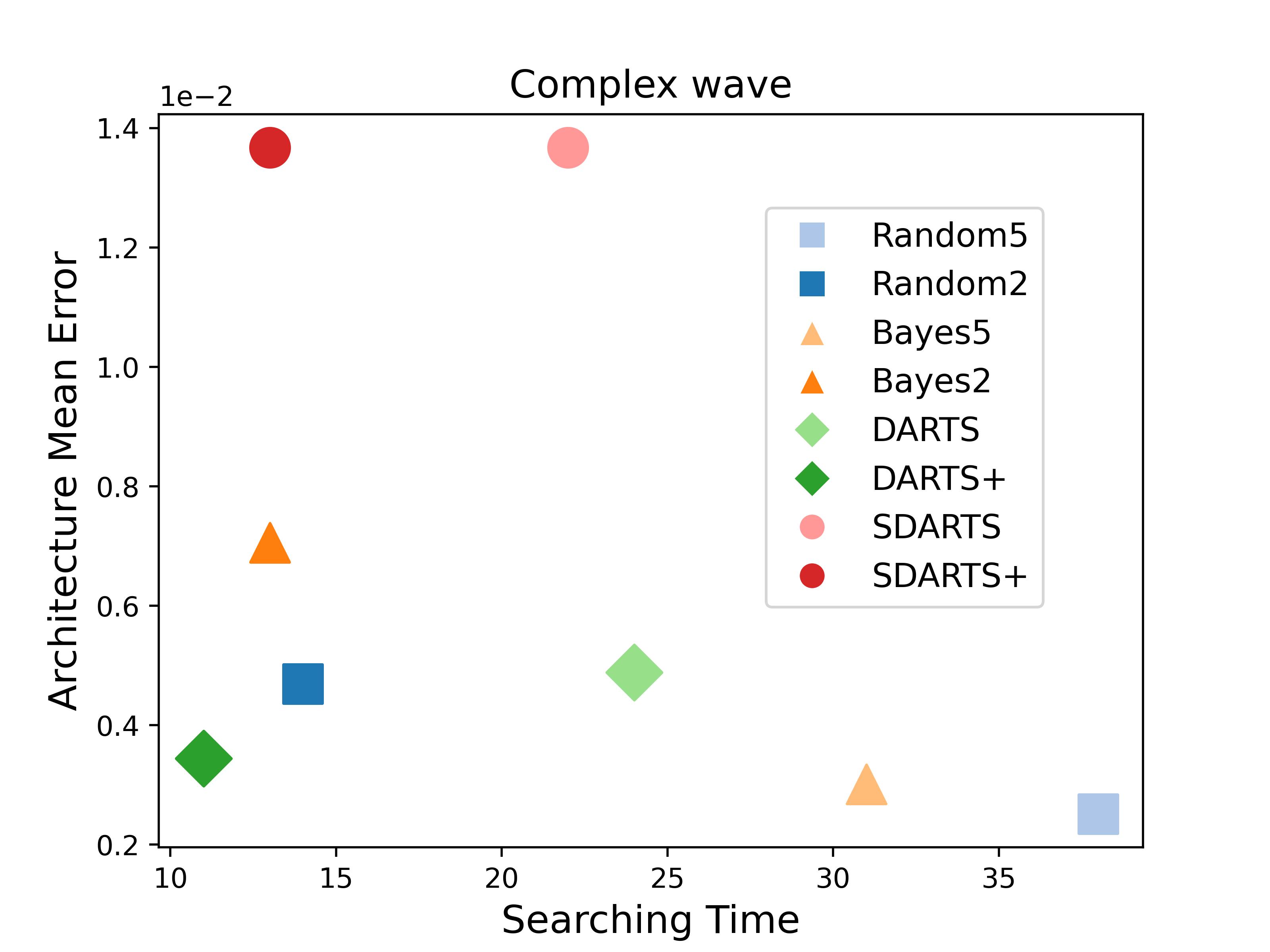}}
    
    \caption{Comparison of PINN-DARTS methods and traditional methods for wave equation with simple and complex solutions.}
    \label{fig4.3}
\end{figure}

For the simple solution case in Figure \ref{sub4.3.1}, PINN-DARTS+ performs better than other traditional search methods with a short search time. PINN-DARTS is slightly worse than traditional search methods. Results for the complex solution case remain unchanged. The optimal network architecture with even widths of simple and complex solutions are similar with details given in Appendix \ref{app1}. Again, we observe that PINN-DARTS and PINN-DARTS+ perform better for problems where the optimal network architecture is deeper and wider.

\subsection{Burgers' Equation}
\indent \indent
Consider the Burgers' equation in \eqref{4.4} with the reference solution obtained by the Chebyshev spectral method \cite{cheby}. 
\begin{equation}
    \label{4.4}
    \left\{
    \begin{aligned}
        & \frac{\partial u(x,t)}{\partial t} + u(x,t) \frac{\partial u(x,t)}{\partial x} 
        - \frac{0.01}{\pi} \frac{\partial^2 u(x,t)}{\partial x^2} = 0 \quad x \in [-1,1],  \ t \in [0,1],  \\
        & u(0,x) = -\sin(\pi x),  \\
        & u(-1,t) = u(1,t) = 0.
    \end{aligned}
    \right. 
\end{equation}
\par 
The training set has 10,000 randomly sampled points inside the domain and 200 randomly sampled points on the boundary. The test set has $256 \times 100$ uniformly sampled points. For the searching phase,  the learning rates of  architecture parameters 
and network parameters are $2 \times 10^{-2}$ and $1 \times 10^{-4}$, respectively. Other parameters are unchanged. Results are presented in Table \ref{tab4.7}.

\begin{table}[htbp]
    \setlength{\belowcaptionskip}{0.1cm}
    \caption{Optimal architecture results of the Burgers' equation. 
    Relative $L^2$ error, search time, and error ratio are averaged over results with 5 initializations.
    }
    \label{tab4.7}
    \centering
    \resizebox{0.7\textwidth}{!}{
    \begin{tabular}{cccc}
        \toprule
        Search Method & Relative $L^2$ Error & Search Time (min) & Error Ratio\\
        \midrule
        Grid Search  & 1.333e-03 &212 & 100.0\%\\
        Random Search, Iter=5  & 2.434e-03 &34 & 182.6\%  \\
        Random Search, Iter=2  & 1.638e-02 &14 & 1228.5\%  \\
        Bayesian Search, Iter=5   & 2.580e-03 &33 & 193.5\% \\
        Bayesian Search, Iter=2   & 7.300e-02 &12 & 5475.1\% \\
        PINN-DARTS   & 2.990e-03 & 31 & 224.3\% \\
        PINN-DARTS+    & 2.001e-03 &14 & 150.1\%  \\
        PINN-SDARTS   & 2.654e-02 & 27 & 1991.0\% \\
        PINN-SDARTS+   & 3.189e-03 & 16 & 239.2\% \\
        \midrule
    \end{tabular}}
\end{table}

\par
\begin{figure}[htbp]
    \centering
    \includegraphics[width=0.55\linewidth]{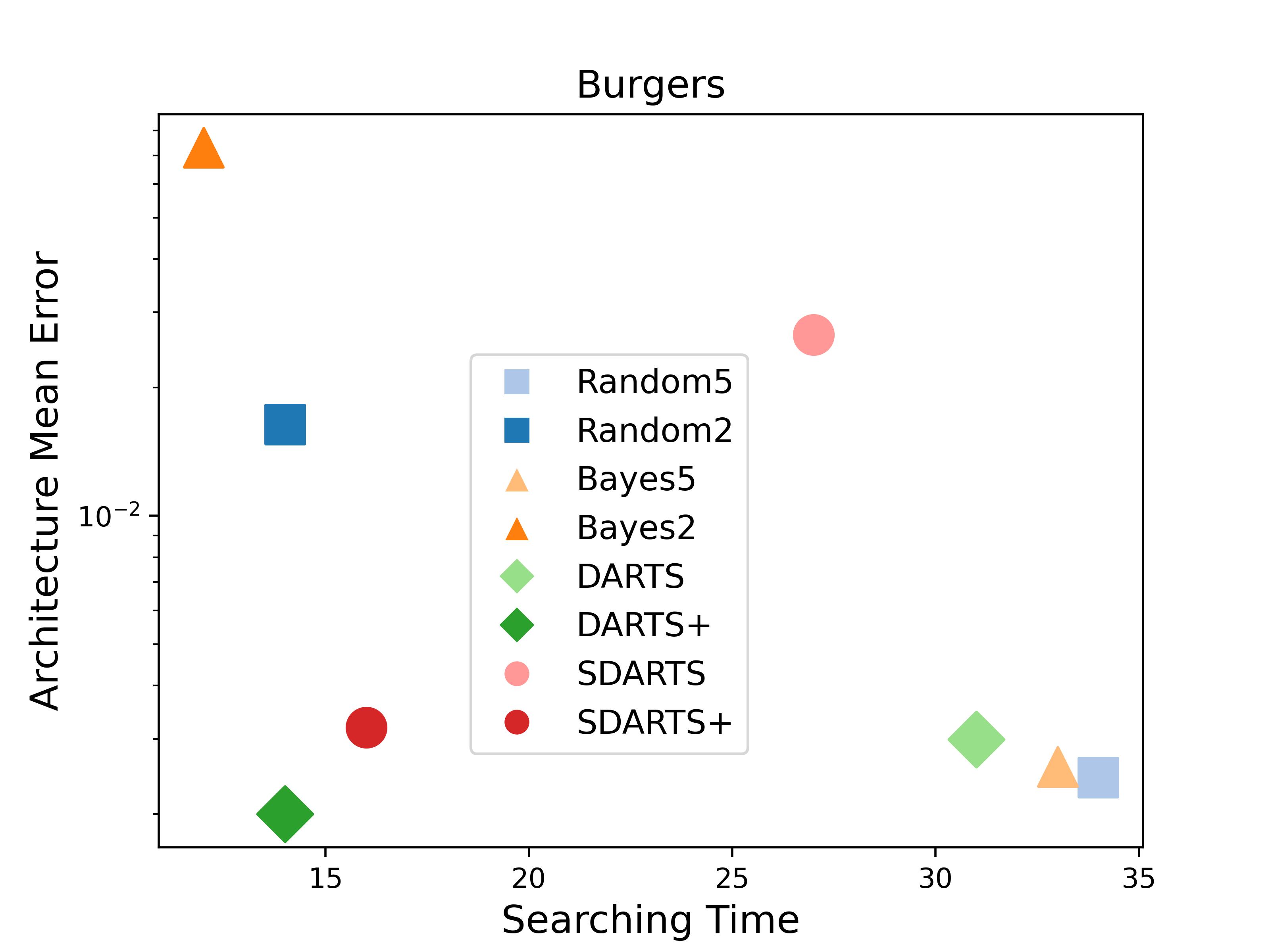}
    \caption{Comparison of PINN-DARTS methods and traditional methods for the Burgers' equation.}
    \label{fig4.4}
\end{figure}

As illustrated in Figure \ref{fig4.4}, PINN-DARTS+ has the smallest average architecture errors. PINN-DARTS+ significantly outperforms traditional search methods within a short search time. 
PINN-DARTS+ shows a higher superiority for the Burgers' equation compared to linear PDEs, suggesting that PINN-DARTS+ performs better on nonlinear PDEs. 
\par
From all experiments, we observe that PINN-DARTS+ performs the best. The detailed architectures found by PINN-DARTS+ and the corresponding numerical solutions are provided in Appendix \ref{app2}.

\section{Conclusions}
\label{section5}
\indent \indent
This paper proposes PINN-DARTS, a class of unsupervised search methods for the optimal network architecture in the framework of physics-informed neural network (PINN), incorporating the widely used differentiable architecture search (DARTS). Both the number of layers and the number of neurons in each layer are allowed to change. The searching phase is unsupervised, only depending on the loss function of PDE residuals without any prior usage of solutions. 
For a variety of PDEs with simple and complex solutions, we find that PINN-DARTS+, which integrates early stopping with PINN-DARTS, has the best performance in terms of search time and architecture accuracy.
The architecture error of PINN-DARTS+ is superior to traditional methods in most numerical experiments with the smallest search time. 
\par
Grid search results indicate that the optimal architectures with even widths differ obviously among different PDEs and solutions. Our study suggests that problem-tailored neural architecture is essential for better accuracy. We also observe that when the optimal network architecture of the current problem has a larger network size, PINN-DARTS+ demonstrates more obvious advantage over traditional search methods, 
possibly due to the superior feature representation ability of architectures with uneven widths.

\par
In future work, given the outstanding performance of PINN-DARTS+ for Burgers' equation, we expect that PINN-DARTS methods will work for more general nonlinear PDEs. This point will be explored more systematically 
from both the numerical perspective and the perspective of approximation theory.

\section*{Acknowledgments}
This work is partially supported by the National Key R\&D Program of China (No. 2022YFA1005200, No. 2022YFA1005202, and No.2022YFA1005203), the NSFC Major Research Plan - Interpretable and General Purpose Next-generation Artificial Intelligence (Nos. 92270001 and 92370205), NSFC grant 12425113, and Key Laboratory of the Ministry of Education for Mathematical Foundations and Applications of Digital Technology, University of Science and Technology of China.
    
\bibliographystyle{elsarticle-num}
\bibliography{ref. bib}

\newpage
\begin{appendices}
\section{Optimal Network Architecture with Even Widths}
\label{app1}
\begin{figure}[htbp]
    \centering
    \subfigure
    {\includegraphics[width=0.49\linewidth]{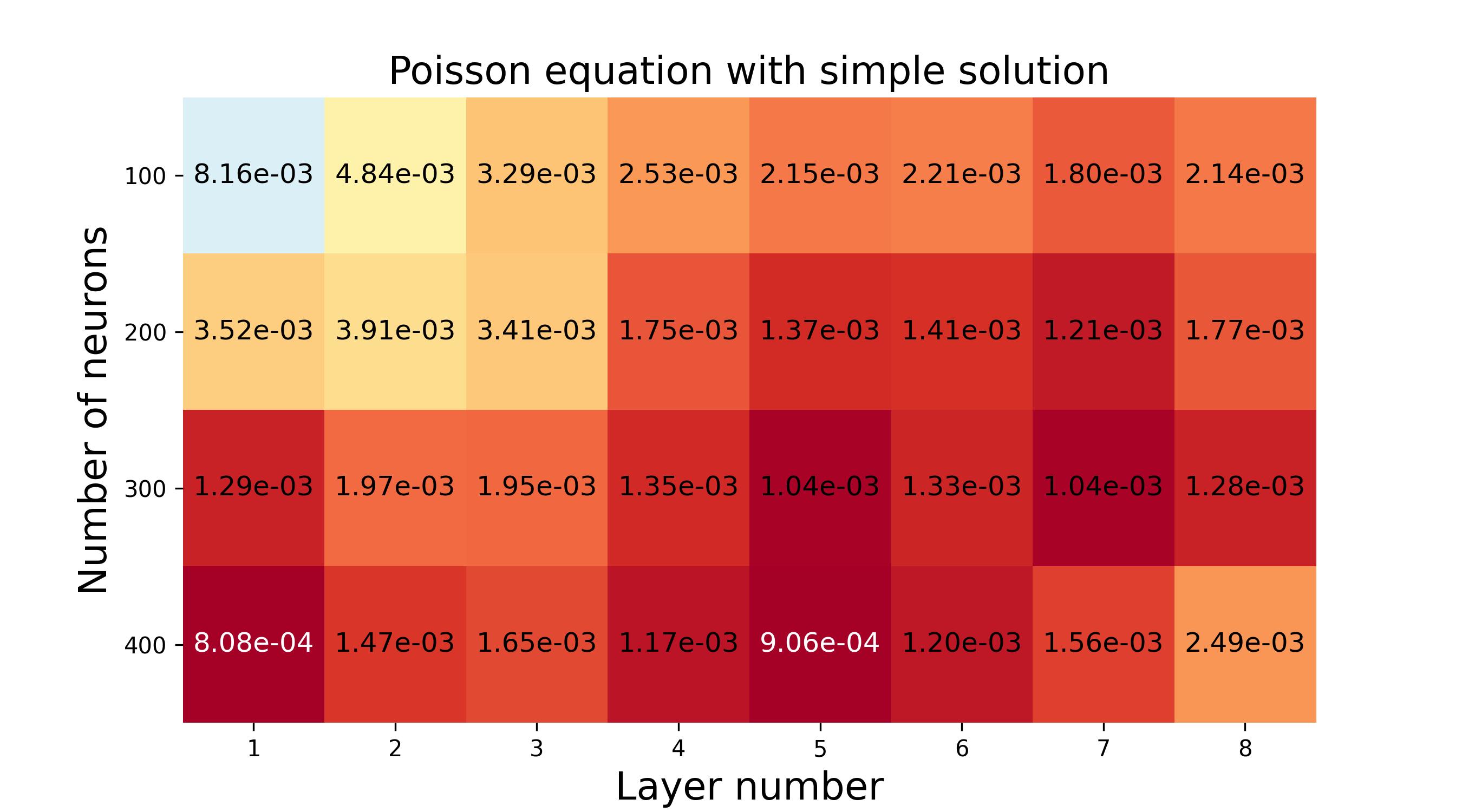}}
    \subfigure
    {\includegraphics[width=0.49\linewidth]{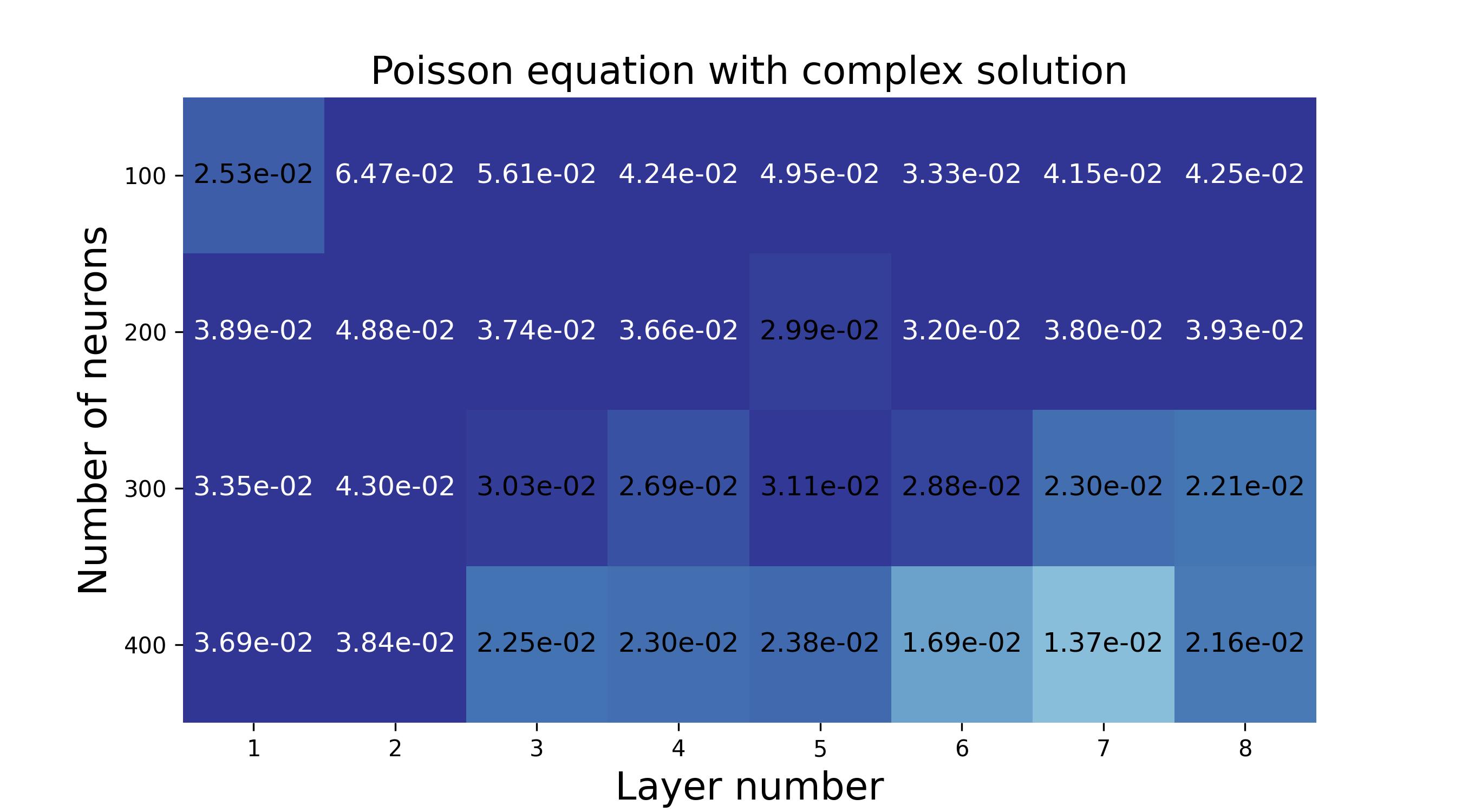}}\\
    \subfigure
    {\includegraphics[width=0.49\linewidth]{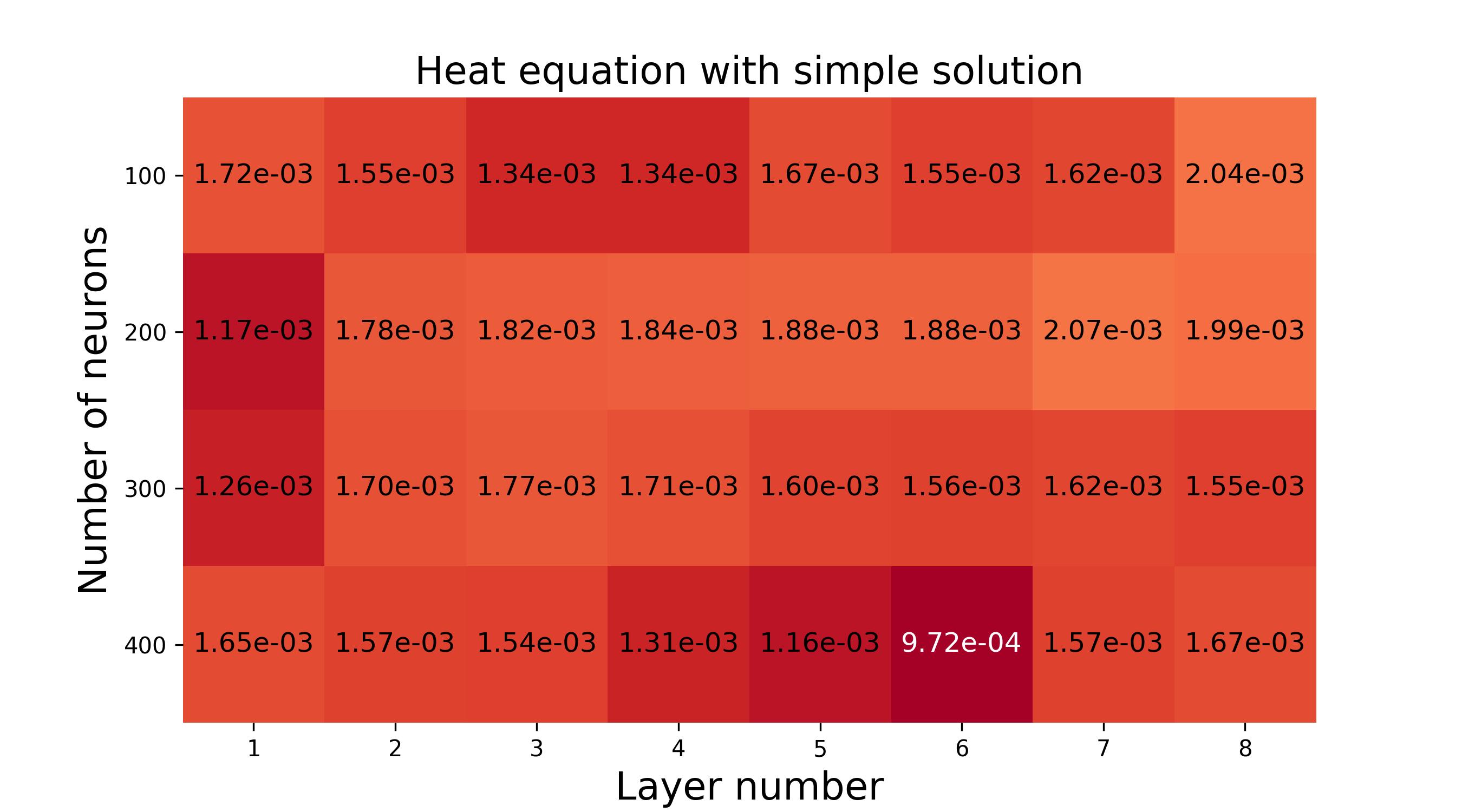}}
    \subfigure
    {\includegraphics[width=0.49\linewidth]{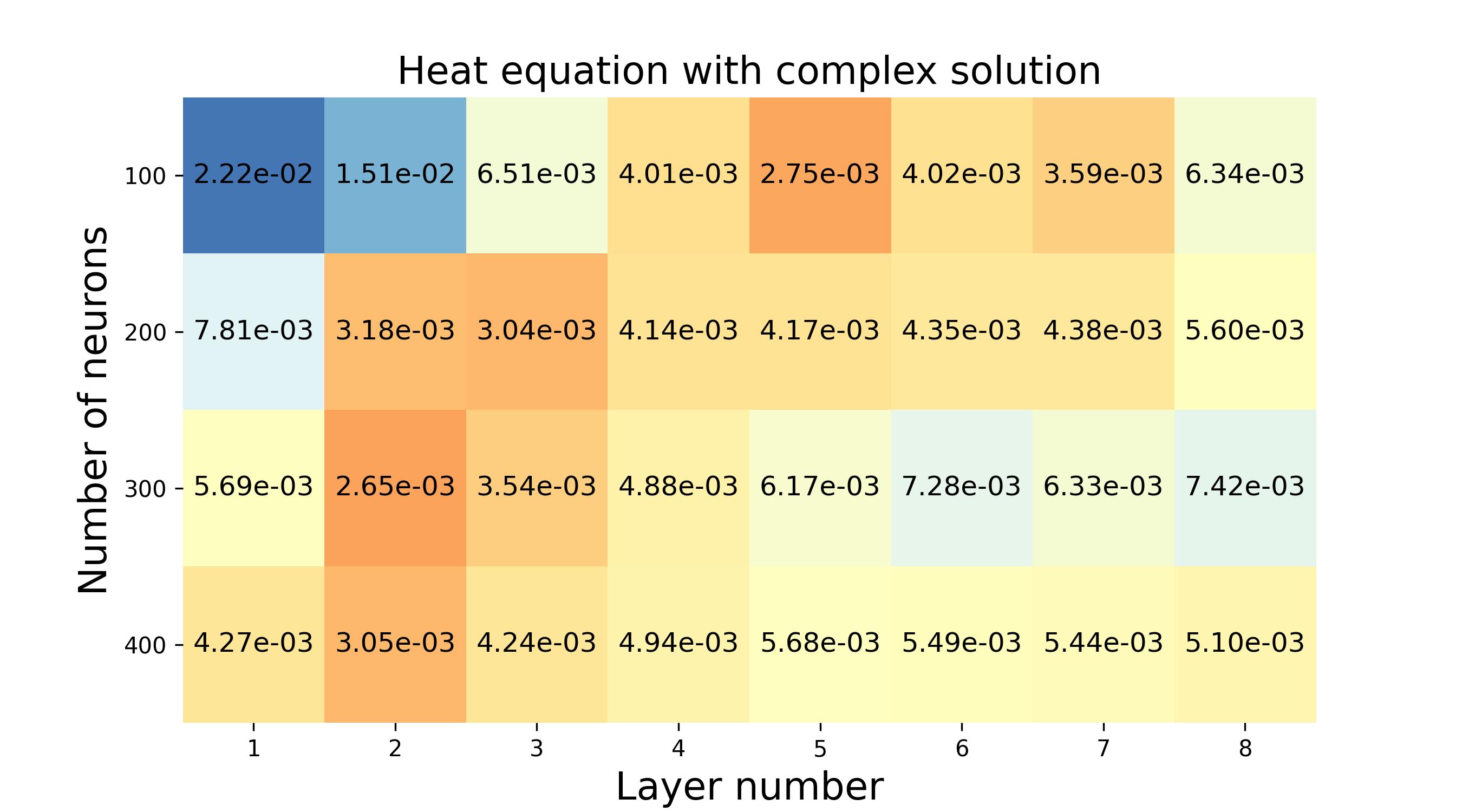}}\\
    \subfigure
    {\includegraphics[width=0.49\linewidth]{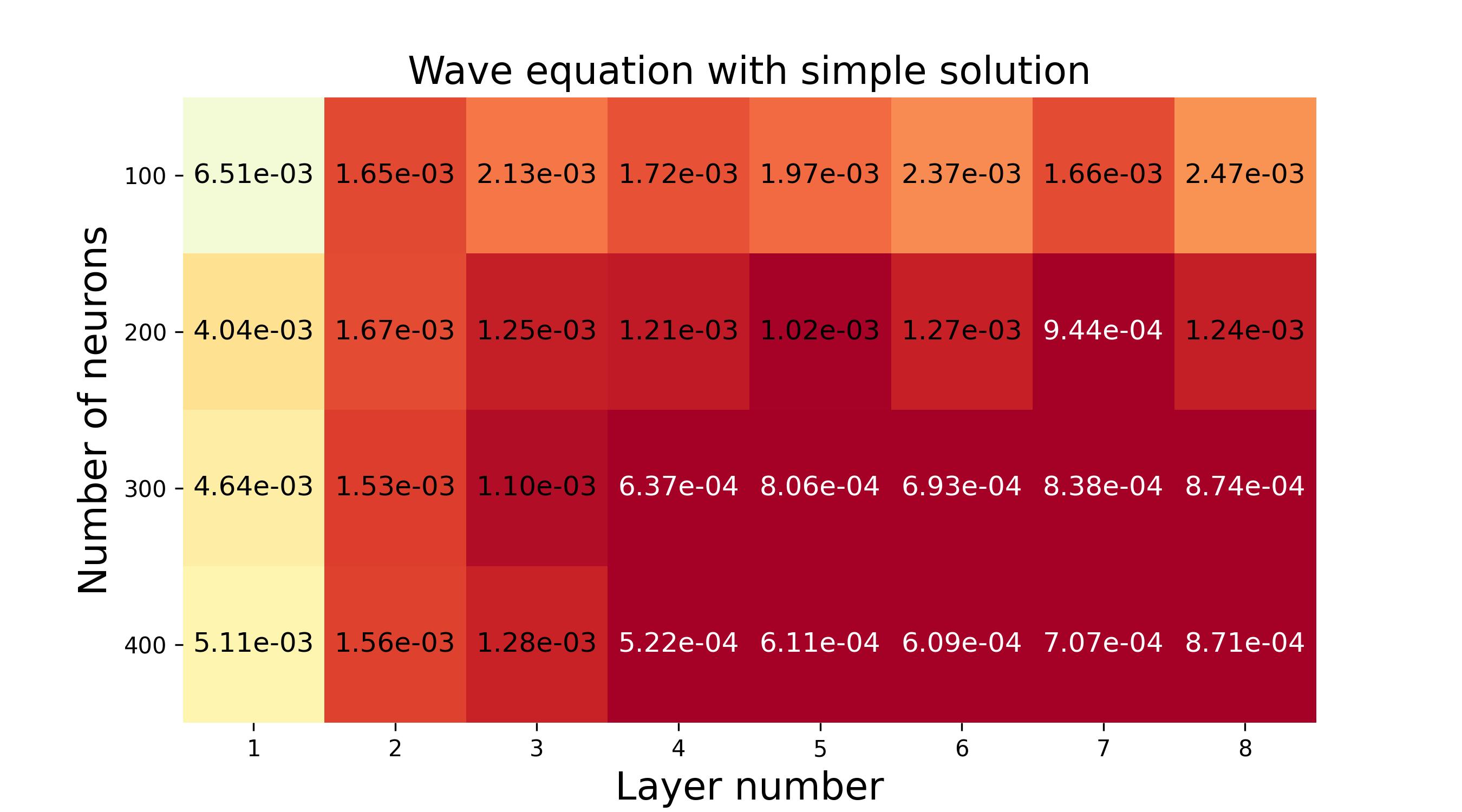}}
    \subfigure
    {\includegraphics[width=0.49\linewidth]{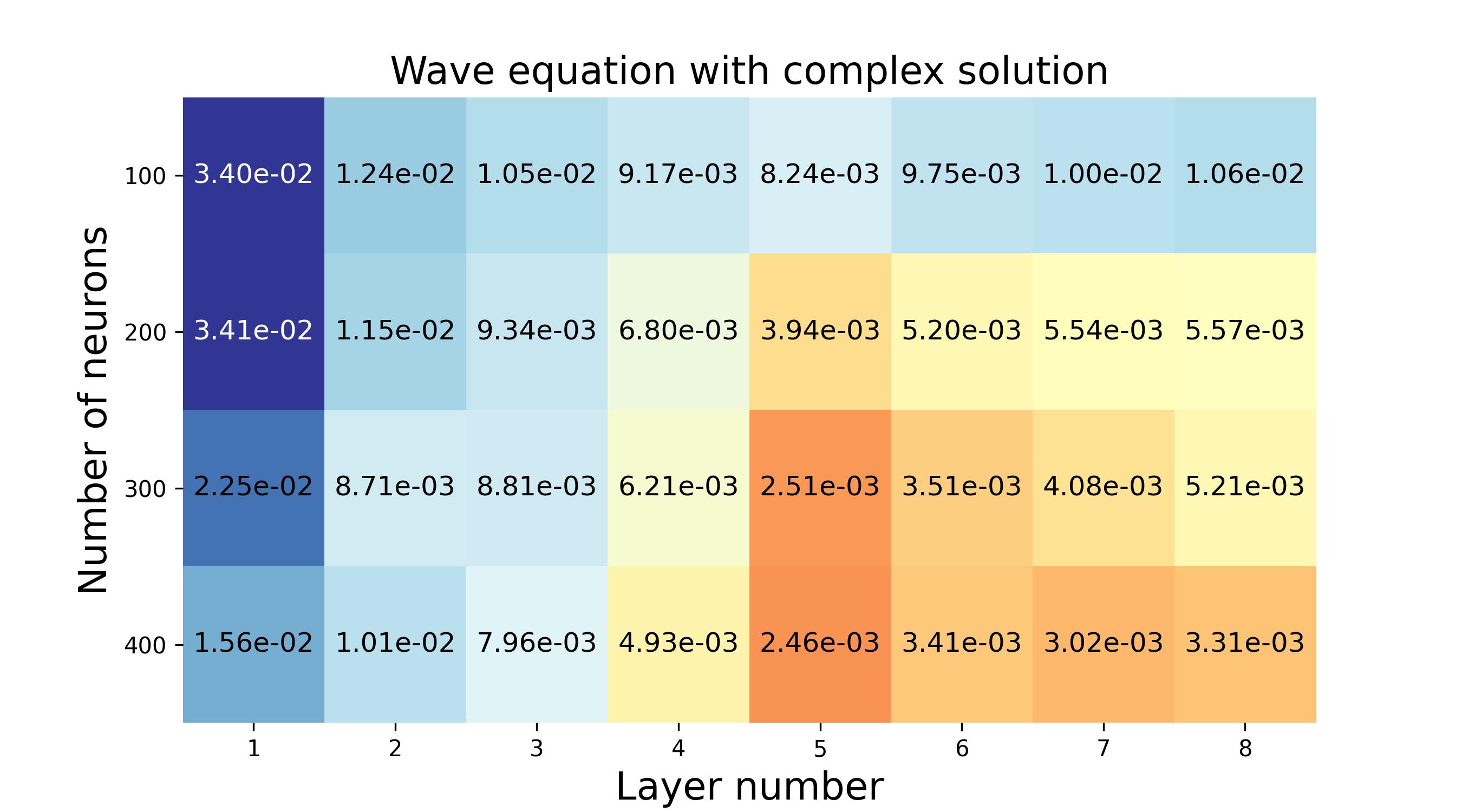}}
    \subfigure
    {\includegraphics[width=0.49\linewidth]{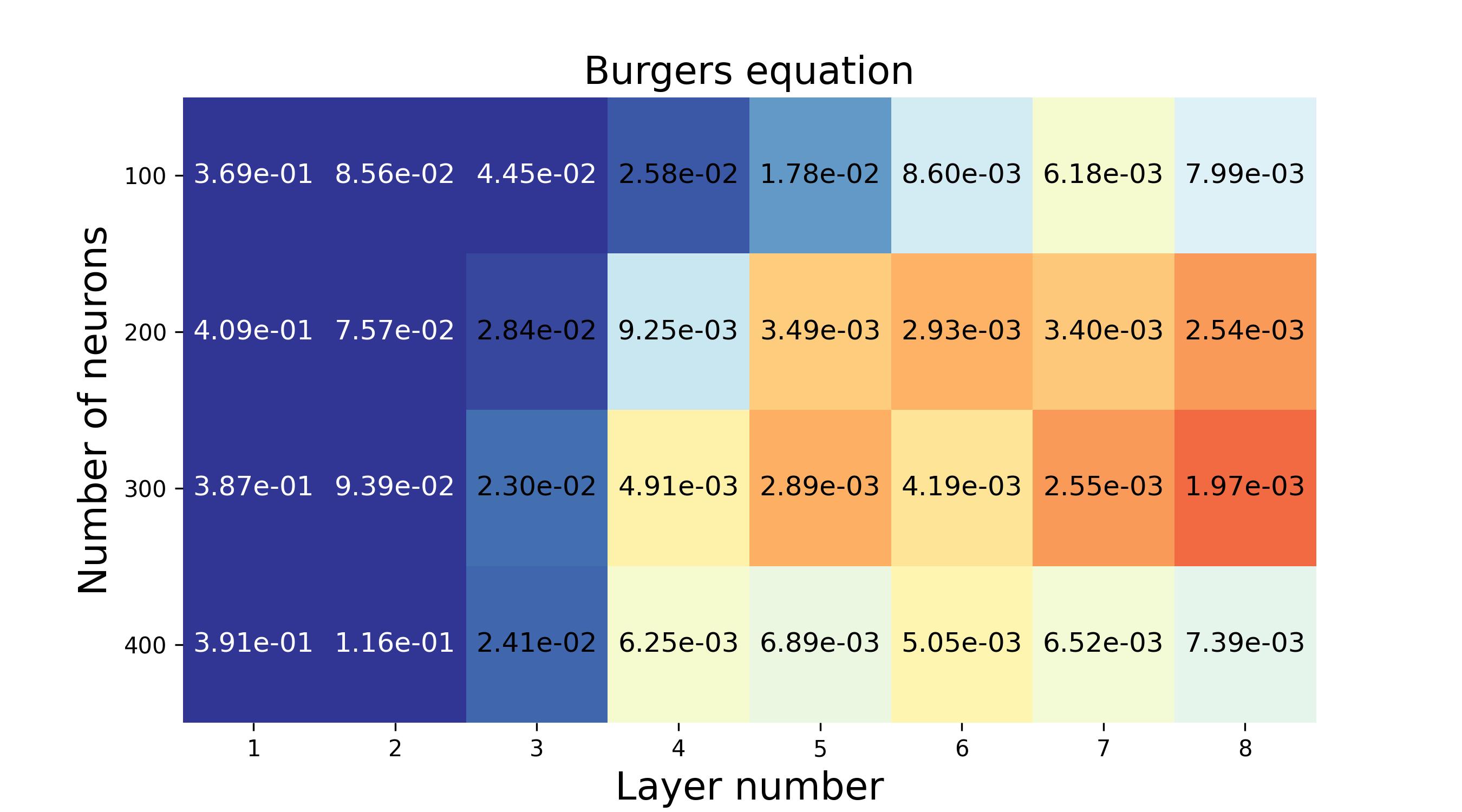}}
    
    \caption{Relative $L^2$ error heatmaps of architectures with even widths. The value in each heatmap cell is the average error of 5 initializations.  Each heatmap corresponds to one PDE with one exact solution. The x-axis is the number of network layers.  The y-axis is the number of neurons in each layer. }
    \label{figa.1}
\end{figure}

\indent \indent
In this appendix,  we present the relative $L^2$ error of all architectures with even widths and the results of the grid search.
Grid search ensures the discovery of the global optimal network architecture.
In Figure \ref{figa.1},  we show the heatmaps of all architectures in different experiments. From the optimal architectures in Table \ref{taba.1},  we find that the optimal architecture of FNN in PINN varies significantly across different PDEs and different solution complexities. The optimal network architecture is tightly related to the PDE and the solution complexity. 

\par
Besides, the grid search results change with different initializations in Table \ref{taba.1}. 
That means the optimal architecture of FNN is also related to the initialization of network parameters. In conclusion, the optimal network architecture is influenced by many factors. 
Problem-tailored network architectures shall be used for better performance rather than any prior knowledge or assumption of the network architecture.

\begin{table}[htbp]
    \setlength{\belowcaptionskip}{0.2cm}
    \caption{Optimal network architecture with even widths of 5 initializations.  Init=Avg corresponds to the architecture with the smallest error, corresponding to the optimal architecture in Figure \ref{figa.1}. We use Width and Depth to characterize these architectures. Equation represents the PDE with different solution complexity.}
    \label{taba.1}
    \centering
    \resizebox{0.77\textwidth}{!}{
    \begin{tabular}{cccc|ccc}
        \toprule
        Equation & Initialization & Width & Depth & Initialization & Width & Depth\\
        \midrule
        Simple Poisson & Avg &  400 & 1
        & Init1 &  400 &  5 \\
        & Init2 & 400  & 1 
        & Init3 & 400  & 6 \\
        & Init4 &  400 & 1 
        & Init5 &  300 & 5 \\
        \midrule
        Complex Poisson & Avg &  400 & 7
        & Init1 & 400 & 6 \\
        & Init2 &  100 & 1 
        & Init3 & 100  & 1 \\
        & Init4 & 400  & 7 
        & Init5 &  100 &  1 \\
        \midrule
        Simple Heat & Avg &  400 & 6
        & Init1 & 400 & 6  \\
        & Init2 & 400  & 8 
        & Init3 & 200  &  1 \\
        & Init4 &  400 &  6 
        & Init5 &  400 &  6 \\
        \midrule
        Complex Heat & Avg &  300 & 2
        & Init1 & 300 & 2 \\
        & Init2 &  100 & 5 
        & Init3 &  100 &  7 \\
        & Init4 & 300  & 3 
        & Init5 &  200 & 6  \\
        \midrule
        Simple Wave & Avg &  400 & 4
        & Init1 & 400 & 4 \\
        & Init2 & 400  & 5 
        & Init3 & 400  & 6 \\
        & Init4 &  400 & 4 
        & Init5 &  400 & 4  \\
        \midrule
        Complex Wave & Avg & 400  & 5
        & Init1 & 300 & 5 \\
        & Init2 & 400  & 5 
        & Init3 & 400  &  6 \\
        & Init4 & 400  & 6 
        & Init5 & 400  & 7  \\
        \midrule
        Burgers & Avg & 300  & 8
        & Init1 & 200 & 6 \\
        & Init2 & 300  & 5 
        & Init3 & 200  &  8 \\
        & Init4 &  400 & 8 
        & Init5 &  300 &  6 \\
        \bottomrule
    \end{tabular}}
\end{table}

\newpage
\section{Search Results of PINN-DARTS+}
\label{app2}
\indent \indent
Among all PINN-DARTS methods, PINN-DARTS+ is the best in terms of search time and architecture accuracy. PINN-DARTS runs much longer than PINN-DARTS+.
The architectural errors of PINN-SDARTS and PINN-SDARTS+ are the largest among all methods. PINN-DARTS+ performs better than traditional methods in most cases. We provide the detailed results of PINN-DARTS+ in this appendix, including the searched architectures and the  corresponding numerical solutions. 
\par
\begin{figure}[htbp]
    \centering
    \subfigure[Simple Poisson ]
    {\label{subb.1}
    \includegraphics[width=0.99\linewidth]{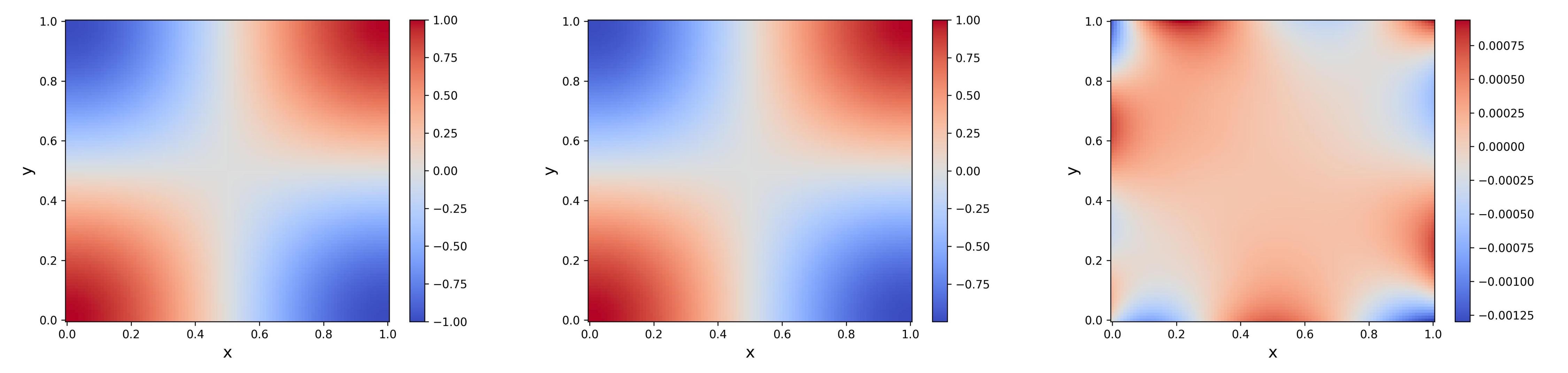}}
    \subfigure[Complex Poisson ]
    {\label{subb.2}
    \includegraphics[width=0.99\linewidth]{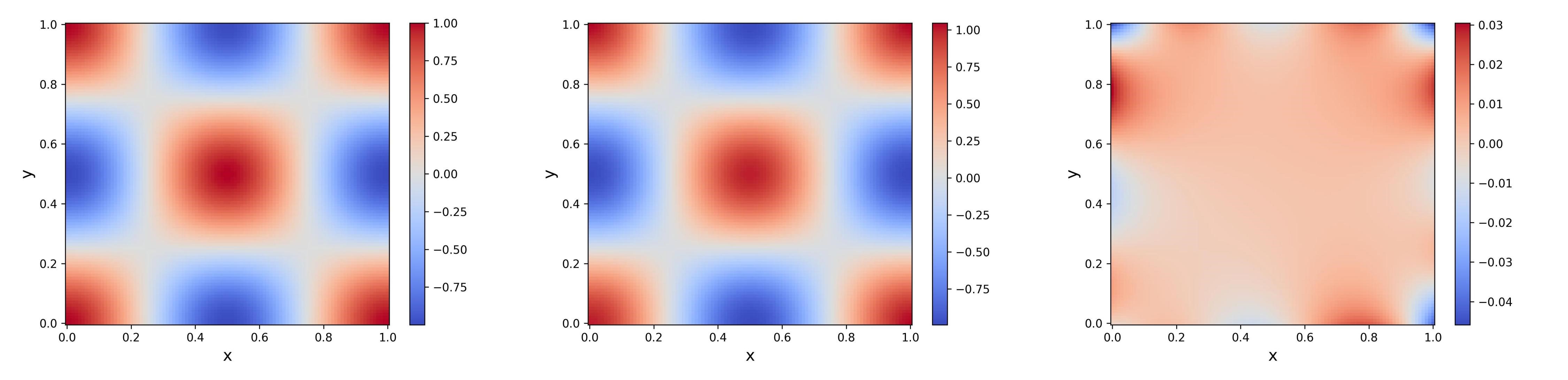}}
    \subfigure[Simple Heat ]
    {\label{subb.3}
    \includegraphics[width=0.99\linewidth]{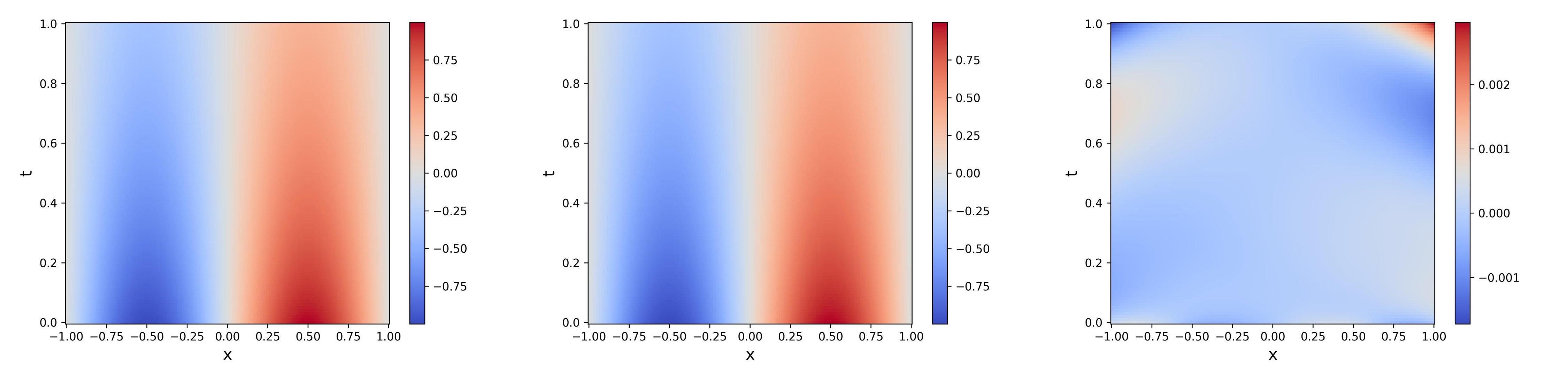}}
    \subfigure[Complex Heat ]
    {\label{subb.4}
    \includegraphics[width=0.99\linewidth]{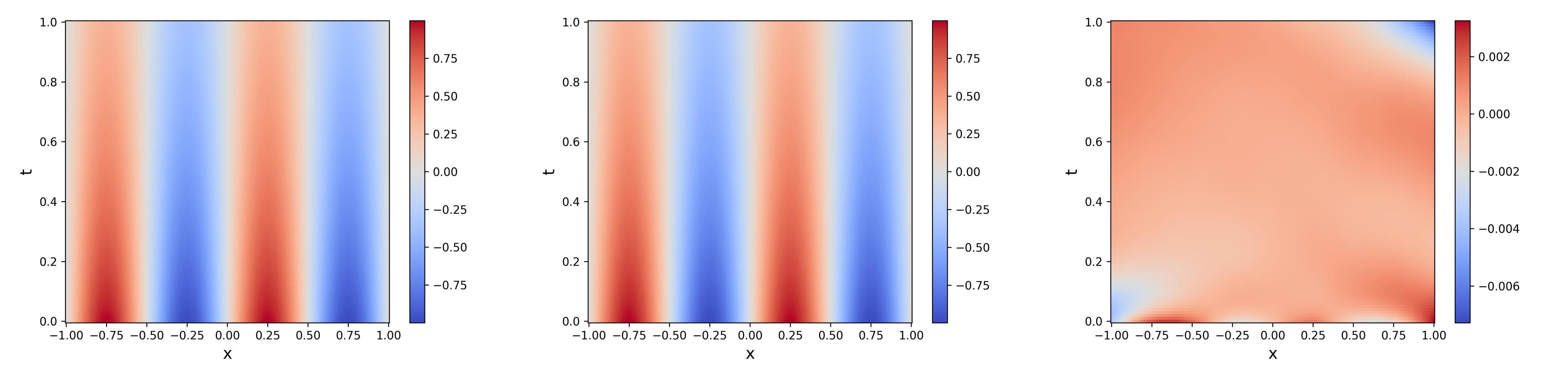}}    
\end{figure}
\begin{figure}[htbp]
    \centering
    \subfigure[Simple Wave ]
    {\label{subb.5}
    \includegraphics[width=0.99\linewidth]{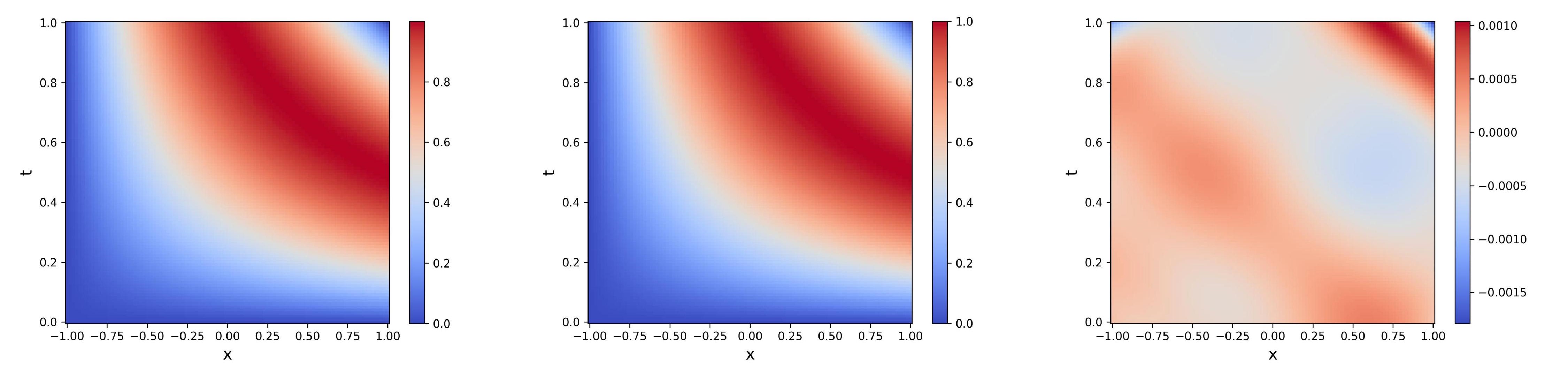}}
    \subfigure[Complex Wave ]
    {\label{subb.6}
    \includegraphics[width=0.99\linewidth]{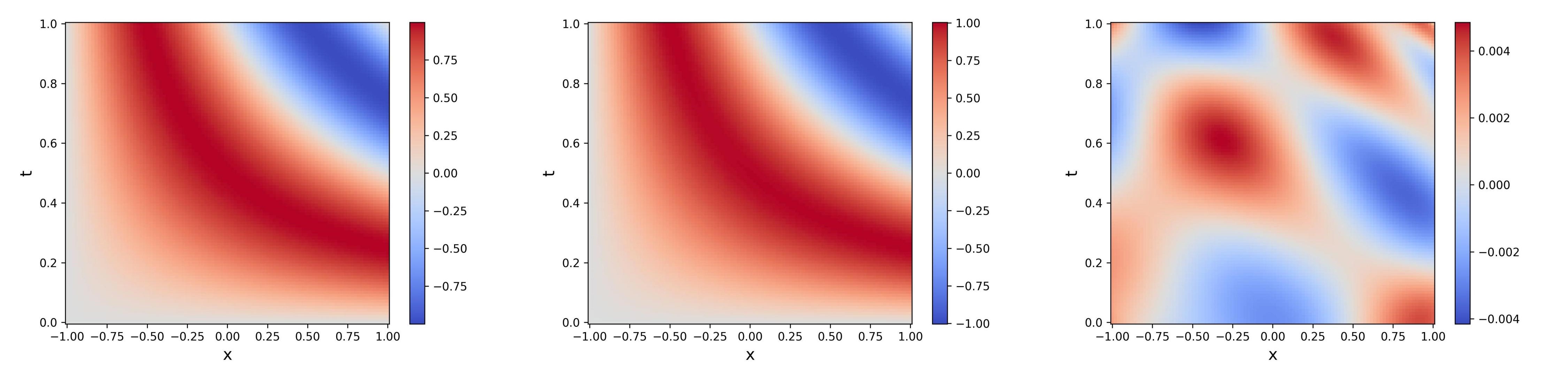}}
    \subfigure[Burgers ]
    {\label{subb.7}
    \includegraphics[width=0.99\linewidth]{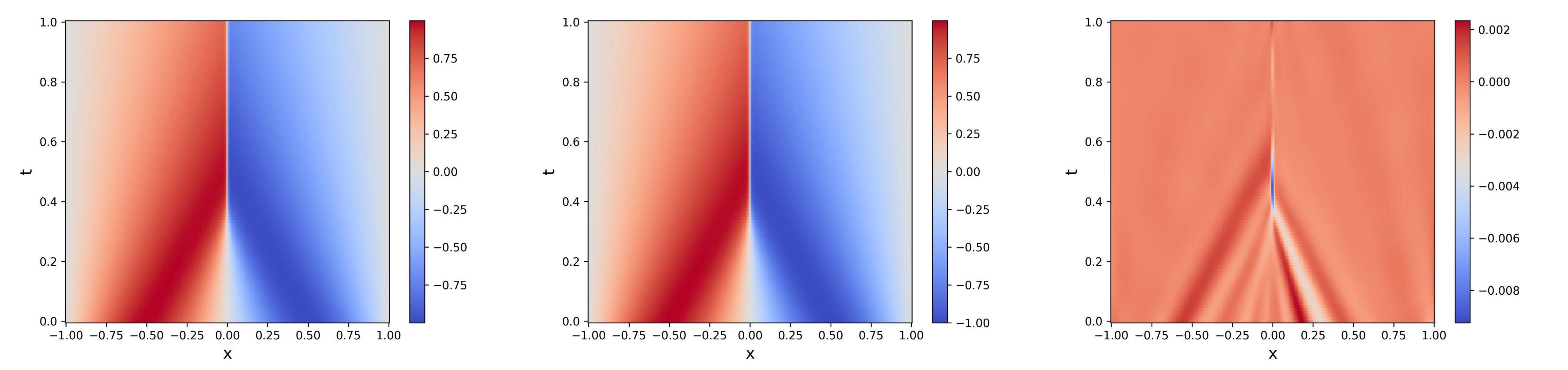}}
    
    \caption{Numerical solutions and error distributions of PINN-DARTS+.}
    \label{figa.2}
\end{figure}

Table \ref{taba.2} presents the search results of PINN-DARTS+ for all experiments, including the searched architectures and loss and error of corresponding numerical solutions. We choose the best searched architectures among 5 initializations for each problem. Architecture is the short form of network hyperparameters. For instance,  (4, 2, 3, 4, 4, 4) represents the network with 6 layers, and the number of neurons for each layer is 400, 200, 300, 400, 400, 400.
Numerical solutions and the corresponding errors are plotted in Figure \ref{figa.2}. Each subfigure corresponds to the searched architecture with the smallest relative $L^2$ error. For example, Figure \subref{subb.1} is the result of Init4 for the Poisson equation with the simple solution. For each subfigure, the left is the exact solution, the middle is the numerical solution of searched architecture, and the right is the distribution of error between the exact solution and the numerical solution. 
\par
Although PINN-DARTS+ performs better than traditional search methods, its architectural accuracy varies with respect to the random initialization of neural network parameters. 
The error of searched architecture under different parameter initializations exhibits large variation. This observation is common in most search methods studied here. 

\begin{table}[htbp]
    \setlength{\belowcaptionskip}{0.2cm}
    \caption{Search results of PINN-DARTS+.}
    \label{taba.2}
    \centering
    \resizebox{0.92\textwidth}{!}{
    \begin{tabular}{ccrcc}
        \toprule
        Equation & Initialization & Architecture &Test Loss& Relative $L^2$ Error\\
        \midrule
        Simple Poisson & Init1 & (4,  2,  3,  4,  4,  4) & 2.283e-05 & 2.220e-03 \\
         & Init2& (3,  4,  4,  4,  3,  3,  4) &1.745e-05  & 8.741e-04 \\
         & Init3 & (4,  4,  4,  4,  4,  3,   4) & 1.709e-05 & 1.422e-03 \\
         & Init4& (3,  4,  4,  4,  3,  4,  4) & 2.481e-05 & 5.459e-04 \\
         & Init5 & (4,  3,  3,  3,  3,  3,  4,  4) & 2.343e-05 & 1.881e-03 \\
        \midrule
        Complex Poisson & Init1 & (1, 1, 1) &2.762e-03 & 1.883e-02\\
        & Init2& (4, 4, 4, 3, 4, 4, 3, 4)  &1.391e-03 & 1.375e-02 \\
        & Init3 & (4, 4, 4, 4, 3, 3, 3, 2) &1.219e-03 & 1.297e-02 \\
        & Init4& (3, 4, 4, 4, 4, 4, 4, 3) &1.487e-03 & 1.411e-02\\
        & Init5 & (3, 4, 4, 4, 3, 4, 4, 3) &1.417e-03 & 1.643e-02 \\
        \midrule
        Simple Heat & Init1 & (3,  2,  3,  3,  4,  4,  4) &1.323e-05 & 1.710e-03 \\
        & Init2& (4,  4,  3,  4,  4,  2,  4) &3.028e-06 & 7.350e-04\\
        & Init3 & (4,  4,  3,  3,  3,  3,  4) &3.987e-06 & 1.223e-03 \\
        & Init4& (4,  3,  4,  3,  4,  4,  4) &4.774e-06 & 9.828e-04\\
        & Init5 & (4,  4,  3,  4,  3,  3) &2.864e-06 & 8.235e-04 \\
        \midrule
        Complex Heat & Init1 & (2,  3,  1,  2) &2.357e-04 & 1.708e-03 \\
        & Init2& (4,  3,  4,  4,  4,  3,  3,  3) &1.835e-04 & 7.197e-03 \\
        & Init3 & (3,  4,  4,  4,  4,  4,  4,  3) &1.268e-04 & 5.482e-03 \\
        & Init4& (1,  1,  1,  1,  1,  1,  1,  3) &7.778e-04 & 3.351e-03 \\
        & Init5 & (3,  4,  4) &1.011e-04 & 4.656e-03 \\
        \midrule
        Simple Wave & Init1 & (2,  2,  3,  2,  3,  3,  4,  4) &4.272e-06 & 7.524e-04 \\
        & Init2& (3,  2) &1.423e-05  & 9.106e-04 \\
        & Init3 & (2,  1,  1,  3,  2,  1,  2,  2) &5.684e-06 & 1.109e-03 \\
        & Init4& (1,  1,  1,  1,  1,  1,  1,  1) &6.412e-06 & 1.337e-03 \\
        & Init5 & (4,  3,  2,  3,  3) &4.654e-06 &  4.745e-04 \\
        \midrule
        Complex Wave & Init1 & (2,  3,  3,  4,  4,  4,  4,  4) &1.018e-04 & 4.959e-03 \\
        & Init2& (4,  4,  4,  4,  4,  4,  4,  4) &8.964e-05 & 3.031e-03 \\
        & Init3 & (4,  4,  4,  4,  4,  4,  4,  4) &8.705e-05 & 3.372e-03 \\
        & Init4 & (3,  3,  3,  4,  3,  4,  4) &8.629e-05 & 2.801e-03 \\
        & Init5 & (4,  3,  4,  4,  4,  4,  4,  4) &8.785e-05 & 3.058e-03 \\
        \midrule
        Burgers & Init1 & (2,  1,  1,  3,  2,  4,  4,  2) &2.919e-05 & 4.055e-03 \\
        & Init2& (4,  3,  1,  4,  1,  3,  4,  4) &1.908e-05 & 3.218e-03 \\
        & Init3 & (3,  3,  3,  3,  3,  2,  2,  2) &6.022e-05 & 1.084e-03 \\
        & Init4 & (3,  3,  4,  4,  2,  4,  4) &7.386e-05 & 9.527e-04 \\
        & Init5 & (4,  3,  4,  4,  3,  2) &6.567e-05 & 6.934e-04 \\

        \bottomrule
    \end{tabular}}

\end{table}


\end{appendices}

\end{document}